\documentclass[3p,review]{elsarticle}
\usepackage{graphicx}
\usepackage{amsmath,amssymb}
\usepackage{array}
\usepackage{hyperref}
\usepackage{color}
\newcommand{\I}{\mathrm{i}}

\newcommand{\sphere}{{\mathbb{S}^2}}

\newcommand{\bx}{{\bf x}}
\newcommand{\by}{{\bf y}}

\newcommand{\bz}{{\bf z}}
\newcommand{\bs}{{\bf s}}
\newcommand{\bd}{{\bf d}}

\newcommand{\bD}{{\bf D}}

\newcommand{\numcalc}{\texttt{NumCalc}{}}

\begin{document}
\begin{frontmatter}{}
\title{Effect of different clustering approaches on the multilevel fast multipole method for the Helmholtz equation }
\author[ARI]{W. Kreuzer\corref{CORRAUTH}}
\ead{wolfgang.kreuzer@oeaw.ac.at}
\author[ARI]{C. Kasess}
\ead{christian.kasess@oeaw.ac.at}

\affiliation[ARI]{
  organization = {Austrian Academy of Sciences, Acoustics Research Institute},
  addressline = {Dominikanerbastei 16},
  postcode = {1010},
  city = {Vienna},
  country = {Austria}
}
\cortext[CORRAUTH]{Corresponding author.}
\begin{abstract}
  The fast multipole method (FMM) is an important component for the boundary element method (BEM), because with the FMM the efficiency and feasibility of the BEM can be enhanced to a large degree. Part of the FMM is grouping the elements of the boundary element mesh into different clusters. The size of these clusters in terms of number of elements and spatial expansion has a huge impact on the efficiency and stability of the method. However, while the theory behind the multipole expansion has been broadly researched, the clustering process itself and its effect on the FMM has been neglected in comparison. Most of the time, for example, it is implicitly assumed that the elements of the mesh have about the same size, which is often not the case in practical applications, e.g., when calculating the sound field around the human head. In this study we compare different types of clustering approaches with respect to stability and efficiency of the underlying FMM applied to meshes that have uniform as well as non-uniform element sizes. Also, some examples are provided for cases where a wrong clustering can lead to numerical problems and instabilities of the FMM-BEM.
\end{abstract}
\journal{Applied Acoustics}
\end{frontmatter}{}

\section{Introduction}
In acoustics the boundary element method (BEM) is an important tool for numerically solving the Helmholtz equation, especially for outdoor scattering problems. Already for moderately sized problems the fast multipole method (FMM) is essential for efficiently using the BEM because with the (single level) FMM the necessary operations the Helmholtz BEM can be reduced from $O(N^3)$ or $O(N^2)$ to $O(N^{3/2})$, where $N$ is the number of unknowns. With the \emph{multi-level} FMM this number can be further reduced to $O(N \log N)$ operations \cite{Coifmanetal93,Rokhlin90,Rokhlin93}.  The memory necessary for computing a BEM solution can be reduced from $O(N^2)$ to $O(N)$ \cite{AmiPro03,CecDar13,Coifmanetal93}.

The main idea behind the FMM is to group neighboring elements of the mesh into different clusters and to reduce the element-to-element interactions to interactions between these clusters. For the single level FMM it can be shown that about $\sqrt{N}$ clusters are quasi optimal \cite{Coifmanetal93}.  For the multi-level FMM the mesh is divided into a cluster tree, where the root consists of one or more clusters. Each of these clusters is again split into at most $2^{\text{Dim}}$ child-clusters that will form the clustering at the next level, where 'Dim' denotes the dimension, in our case Dim = 3. New levels are added to the cluster tree until the number of elements per cluster in the last level is about 25.

While the theoretical background and the splitting of the respective integral kernels involved with the (in this case) Helmholtz BEM have been thoroughly investigated in the past \cite{AmiPro00}, the actual process of clustering and its effect on the FMM has been neglected in comparison. For the multi-level FMM  there is some research on how to build up the cluster tree, see for example  \cite[Chapter 7]{SauSch10} or \cite{Chengetal06,Kreuzeretal24,Lietal24,Yunetal20}. These approaches range from the most common approach of an oct-tree to a grouping based on k-means clustering. In connection with the oct-tree clustering there are some approaches that start with the full mesh at root level, while others already use a clustered mesh at the root that is then further grouped into smaller child clusters. However, in general, only the algorithm that is used to build up the cluster tree is described, but, to our knowledge, there are hardly any comparisons between these approaches in terms of practicality, stability, and efficiency of the fast multipole method connected to the clustering.

Almost all clustering methods assume a mesh where the sizes of the elements are about the same. This is often not the case for practical applications, e.g, if the mesh is refined adaptively or if the mesh consists of two parts that are discretized differently, see for example \cite{Brinkmannetal23, Ziegelwangeretal16}. This will lead to a cluster tree where on a given level the number and/or the radii of the clusters (i.e. the biggest distance of a vertex in the cluster to the cluster midpoint) may differ to some degree. Even if a mesh has elements of about the same size, the number of elements per cluster may vary considerably in certain situations as will be illustrated with the example of a clustering of the unit-sphere in Section \ref{Sec:Problems}.

In this study we will look at practical considerations in connection with different clustering approaches for uniform and non-uniform meshes. For the uniform mesh we will use the standard benchmark problem of scattering of a plane wave from a sound hard sphere. For this case we will look at the effect of the number of clusters in the root level of an oct-tree that are either generated by simply dividing the mesh along the three main coordinates or by using a k-means clustering based on the vertices of the mesh. For a practical benchmark for  non-uniform meshes we will look at the calculation of sound scattering from the human head, where the outer ear (pinna) is discretized with very small elements whereas the rest of the head uses coarser elements. Especially in connection with this example we will also look at the influence of the mixture of large and small elements on the efficiency and the stability of the FMM. We will identify some conditions where the clustering of the mesh may cause numerical problems for the FMM that may occur for non-uniform as well as uniform meshes.

The manuscript is structured as follows. Although already described in several publications we give a short summary of the multipole method in Section \ref{Sec:FMM} to  keep the paper self-contained. In Section \ref{Sec:Clusteringmethods} we will describe details of the clustering  methods used for the numerical experiments described in Section \ref{Sec:Experiments}.   In Section \ref{Sec:Conclusion} we give a short summary and discuss the results of the numerical experiments.

\section{The Fast Multipole Method}\label{Sec:FMM}
The FMM aims at speeding up matrix vector multiplications leading to a massive reduction in computation time. Also, the system matrix does not need to be stored in full leading to an impressive reduction in memory requirements. In a nutshell, instead of calculating interactions between each pair of unknowns necessary for the BEM, these interactions are reduced to groups of unknowns. To that end the mesh is divided into different clusters (i.e. groups of elements close to each other), and interactions between individual unknowns are split into local expansions in each single cluster, and interactions between the different cluster midpoints which are, in general,  defined as the mean coordinates of all vertices in the respective cluster. For each cluster its radius is defined as the largest distance between the midpoint of the cluster and vertices contained in the cluster.
\begin{figure}[!h]
  \begin{center}
  \includegraphics[width=0.5\textwidth]{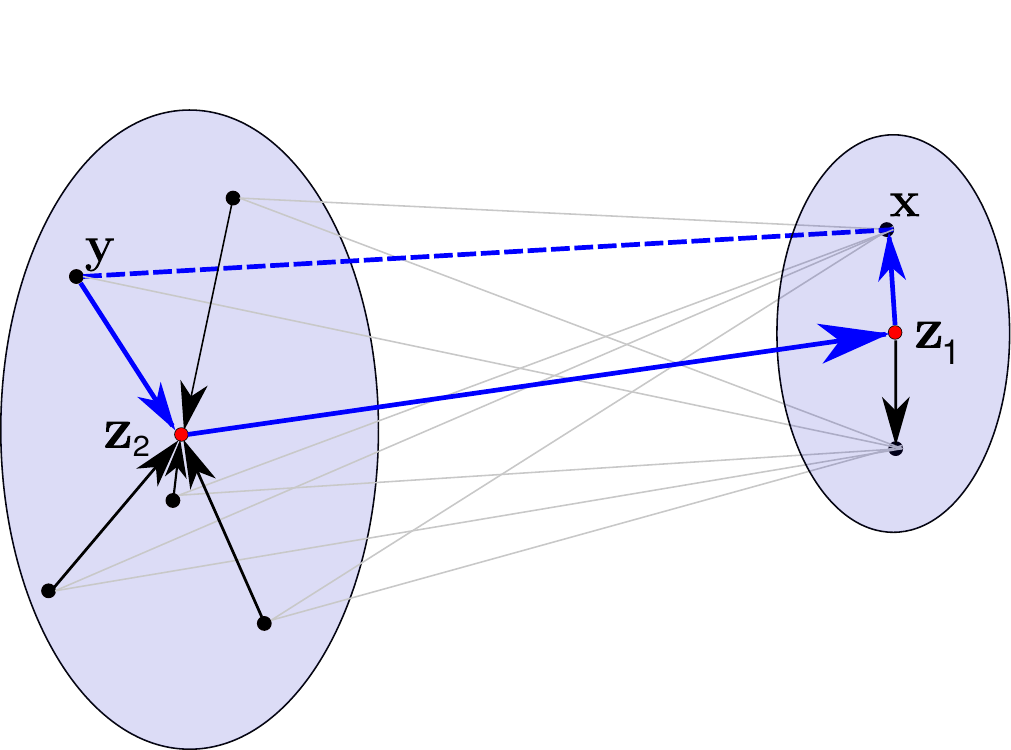}
  \caption{FMM Scheme.  Instead of calculating interactions between all nodes (gray lines), interactions are reduced to their cluster midpoints.} \label{Fig:FMMScheme}
  \end{center}
\end{figure}

The multipole expansion of the Helmholtz Green's function  is based on its representation using spherical harmonics \cite[Eq.~(13)]{Coifmanetal93}:
\begin{align*}
\frac{e^{\I k ||\bx - \by||}}{4\pi ||\bx -  \by||}
&\approx \frac{\I k}{4\pi}
\int\limits_{\sphere}e^{\I k(\bx - \bz_1 + \bz_2 - \by)\cdot \bs}
\sum\limits_{\ell=0}^L h_\ell(k||\bz_1 - \bz_2||)
P_\ell\left(
\frac{(\bz_1 - \bz_2) \cdot \bs}{||\bz_1 - \bz_2||}
\right) d\bs\\
&\approx \frac{\I k}{4\pi}
\int\limits_{\sphere}e^{\I k\bd \cdot \bs}
\sum\limits_{\ell=0}^L h_\ell(k||\bD||)
P_\ell\left(
\frac{\bD \cdot \bs}{||\bD||}
\right) d\bs,
\end{align*}
where $\sphere$ is the surface of the unit-sphere, $h_\ell$ are the spherical Hankel functions of order $\ell$, and $P_\ell$ denotes the Legendre polynomial of order $\ell$. The exponential terms $e^{\I k \bd \cdot \bs} = e^{\I k(\bx - \bz_1) \cdot \bs}e^{\I k(\bz_2 - \by) \cdot \bs}$ can be interpreted as local expansions in two clusters around their respective midpoints $\bz_1$ and $\bz_2$.
The vector $\bD = \bz_1  - \bz_2$ is the ``global'' distance vector between the midpoints of the two clusters. The truncation parameter $L$ of the multipole expansion can be set following the semi-empirical formula \cite{Coifmanetal93}
\begin{equation}\label{Equ:Truncation}
L  = 2kr_{\max} + \gamma \ln\left( 2k r_{\max} + \pi \right),
\end{equation}
where $r_{\max}$ is the biggest radius of the clusters at a given level, and $\gamma$ is a parameter that depends on the desired accuracy for the approximation. The value of $L$ implicitly determines the number of quadrature nodes $\bs_i$  necessary for the numerically calculating the integral over $\sphere$~\cite{Rahola96} and is subject to the following general rules~\cite{Coifmanetal93,Darve00}:
\begin{itemize}
\item $L$ and the number of quadrature nodes should be \emph{high enough} to accurately represent the oscillating local expansions $e^{\I k \bd\cdot \bs}$, see \cite[Theorem 1]{Greengardetal98}.
\item On the other hand, $L$ should be \emph{small} such that $h_\ell(k||\bD||)$, which has a singularity for $k||\bD|| \rightarrow 0$, can be numerically calculated with sufficient accuracy. A high order of the spherical Hankel function implies that $k||\bD||$ needs to be relatively large to ensure numerical stability. This in turns means that clusters for which the  multipole method can be used need to be sufficiently far from each other for the FMM to be numerically stable.
\end{itemize}
The first point can be used as a motivation for the dependence of $L$ on the cluster radius in Eq.~(\ref{Equ:Truncation}), the second point has influence on the criterion to define when the FMM can be applied to a cluster pair (the far-field  criterion \cite{Coifmanetal93}): If for two clusters $\mathcal{C}_1$ and $\mathcal{C}_2$ with midpoints $\bz_1$ and $\bz_2$, and radii $r_1$ and $r_2$
\begin{equation}\label{Equ:FarfieldCrit}
  ||\bz_1 - \bz_2|| > \frac{2}{\sqrt{3}}(r_1 + r_2 ), 
\end{equation}
both clusters are within each others far field and the multipole expansion can be applied. Otherwise, they are within each others near field  and the ``classical'' BEM approach needs to be used to calculate the entries in the system matrix for all elements in this cluster pair.

For the multi-level FMM (MLFMM) the mesh is clustered using a hierarchical tree structure. A common approach to build such a tree is to define the root of the tree as one single cluster containing all elements of the mesh. Based on the coordinates of the vertices contained in the cluster, a box is defined around the cluster, that is then subdivided into $2^\text{Dim}$ sub-boxes, where 'Dim' denotes  the spatial dimension of the problem. The elements in each of these sub-boxes then define the new clusters  at the next level. The original cluster is called parent, the newly generated clusters are its children. This procedure is repeated for all clusters at a given level until the number of elements contained in a cluster is about 25. The clusters at the highest level (i.e. those clusters without children) are commonly defined as leafs of the cluster tree. The advantage of using a cluster tree is that
\begin{itemize}
\item the near field needs to be calculated only for the clusters at the leaf level which, in general, contain few elements,
\item the FMM for a cluster on a given level is only done with those clusters that are in its interaction list, i.e., clusters that are children of the parents neighboring clusters. For the rest of the clusters the local expansion is passed up one level up and the procedure is repeated at the parent level.
\end{itemize}
This results in an $O(N \log N)$ algorithm for a matrix--vector product, where $N$ is the number of unknowns. In general, it can be said that the \emph{number} of elements in the leaf clusters define the size of the near-field matrix, which mainly influences the efficiency, whereas the \emph{radii} of the clusters are related to the expansion length of the FMM at each level, which also has influence on the stability of the method. For more details we refer to \cite{Kreuzeretal24}.
\section{Clustering methods}\label{Sec:Clusteringmethods}
In the following we look at three different approaches for defining the cluster tree. 
\subsection{Box1 Clustering}\label{Sec:Box1}
Variants of this approach are for example used in \cite{Kreuzeretal24, Lietal24}. The main idea is to already partition the mesh at the root level into multiple equally sized clusters. This approach has the advantage that the radii and in turn the expansion length of the FMM can be kept relatively small which, in turn, may enhance stability and efficiency. 

The root of the cluster tree is generated by putting a box around the whole mesh that is then subdivided into equally sized boxes with fixed edge length $h_0$. In the open source software \numcalc~\cite{Kreuzeretal24}, that is the basis for the numerical experiments in this study, this initial edge length can be provided by the user or left at a default value of
\begin{equation}\label{Eq:DefaultBox}
h_0 = \left( \left( \frac{{N}^{1/2}}{0.9} \right)^{1/2} l_0 \right),
\end{equation}
where $l_0$ is the square root of the mean element area.

The motivation behind this value lies in the suggestion for the single level FMM to use about $\sqrt{N}$ clusters, where $N$ is the number of elements \cite{Coifmanetal93} . This choice can be motivated as follows: Assume a mesh with $N$ elements and $n_0$ clusters. The most costly parts of the FMM are the near-field interactions, that can be estimated by the number of the elements in each nearfield cluster pair times the number of clusters: $O\left( \left(N/n_0\right)^2  n_0 \right)$. The cluster-to-cluster interactions depend on the number of clusters $n_0$ and the number of quadrature nodes on the sphere. The effort can in principle be approximated by $O(n_0^2 L^2)$ where $L$ is the expansion length defined by Eq.~(\ref{Equ:Truncation}). The cluster radius scales with $1/\sqrt{n_0}$, the wavenumber can be assumed to scale approximately with $\sqrt{N}$ (see, e.g.,  \cite{AmiPro03}). Under the assumption that $kr \approx \sqrt{N/n_0}$ is big enough that the log term in Eq.~(\ref{Equ:Truncation}) can be neglected the whole effort scales with $O(N^2 / n_0) + O( n_0 N )$. The $O(N^2/ n_0)$ term also contains the estimation for the local expansions. It is then easy to show that $n_0 = O(\sqrt{N})$ minimizes the above expression. For deriving the initial box length for the SLFMM we assume that the mesh is relatively smooth and that it covers an area of $A = N A_0$ where $N$ is the number of elements and $A_0$ is the average element area. Each of the $n_0 = \sqrt{N}$ clusters, should ``cover'' a part of the \emph{surface} mesh with area $A_{\mathcal{C}} = A / n_0 = A_0 \sqrt{N}$, thus it is a good assumption to set the initial box length to $\sqrt{ A_0 \sqrt{N}}$.
In this rough estimate the log term in the calculation of the truncation length $L$ is neglected, which is, in general, only valid for large $kr$.
%This rough estimate is based on the assumptions that the log term can be neglected in the estimation of the initial box length is, in general, only valid for large $kr$. %Fig.~\ref{Fig:SLFMM} shows the computing time and the peak memory for the example of the single level FMM applied to the uniform sphere mesh with $N = 69620$ elements at $f = 1500$ Hz. 
%\begin{figure}[!h]
%  \includegraphics[width=0.46\textwidth]{Pics/WtimeSLFMM}
%  \includegraphics[width=0.46\textwidth]{Pics/MemSLFMM}
%  \caption{Computing time and peak memory for the single level FMM applied to the uniform mesh of the unit sphere at $f = 1500$ Hz.}\label{Fig:SLFMM}
%\end{figure}
%The default box length $\sqrt{ A_0 \sqrt{N}} \approx 0.218$\,m is marked with a circle. For that box length the number of clusters is given by 314, which is slightly larger then $\sqrt{N} = 263$. It can be seen that for this specific example the optimal initial box length is $h_{\text{opt}} \approx 0.15$\,m that corresponds to 689 clusters. In this case the above estimation is not valid, partly because the log factor in the estimation of the truncation length $L \approx kr + \gamma \log(kr + \pi)$ cannot be neglected. 

For the MLFMM a rough estimation for a good box length is more complex. The  different levels imply different expansion lengths and different quadrature nodes on the sphere. In this case the logarithm-term used for the calculation of the expansion length cannot be neglected anymore. An additional  complexity is introduced by the rounding involved in the calculation of the number of necessary levels, see Eq.~(\ref{Equ:Nlevels}). 

Nevertheless,  numerical experience has shown that setting the number of clusters at the root level to about $\sqrt{N}$ is also a good choice for the MLFMM. The 0.9 factor was introduced in Eq.~(\ref{Eq:DefaultBox}) because numerical experience showed that this may enhance the performance slightly, see also the numerical experiments in Section~\ref{Sec:Experiments}.

The actual edge length of the sub-boxes will differ from $h_0$ because the number of subdivisions of the original box is naturally defined by \emph{rounding} the quotient of box edge length and $h_0$, thus the actual boxlength is given by:
\begin{center}
  %\texttt{Edge-length(Box$_0$) $\cdot$  round( h$_0$ / Edge-length(Box$_0$) )}.
  Edge-length(Box$_0$) $\cdot$  round( h$_0$ / Edge-length(Box$_0$) ).
\end{center}
This also implies that, depending on the spatial dimensions of the mesh, the number of subdivisions in the three different coordinate directions may differ from each other. Note that in the following $h_0$ always denotes the user provided or default initial edge length \emph{before} the rounding process.

The elements in each box form one single cluster, however, if a cluster is too small in terms of cluster radius it is merged with its nearest neighbor. The idea behind this step is to avoid numerical problems that could be caused by having very small and large clusters at a level of the FMM tree, see Section \ref{Sec:Problems}. 

For generating the  children, a bounding box is drawn around the elements of the parent cluster, which is then subdivided into equally sized boxes with edge length $h_\ell = 0.5^{\ell} h_0$, where $\ell$ is the current level in the cluster tree. As a consequence, a parent cluster in 3D can have theoretically up to 8 children. If a parent box at level $\ell$ has an edge with a length smaller than $0.5^{\ell+1} h_0$ no subdivision will be made along this edge. This  implies that it is possible for a parent cluster to be its own child. 

For practical reasons, the number of levels is often defined beforehand using the quotient of the average number of elements in the root clusters and the average number of elements in the leaf clusters:
\begin{equation}\label{Equ:Nlevels}
  {\text{Levels}} \approx \text{round}
  \left(
     \log_4
     \left(
     N_{\text{root}} / N_{\text{leaf}}
     \right)
   \right)
 \approx
 \text{round}
 \log_2\left(
h_0 / (5 l_0)
  \right),
\end{equation}
where the last term is based the  initial boxlength estimation used in \numcalc, while the 2nd one is based on number of clusters.

The determination of the cluster levels, construction of the cluster tree, and the stability of the FMM are based on the assumptions that
\begin{itemize}
\item the geometry of the scatterer is relatively smooth,
\item all elements have about the same size,
\item all clusters on a single level contain about the same number of elements,
\item all clusters on a single level have about the same radius,
\item the extensions of each child cluster in the $x,y,z$-directions are about half the extensions of its parent cluster.
\end{itemize}
\subsection{Box2 Clustering}
This recursive method uses a slightly different approach than the Box1 described in Section \ref{Sec:Box1}. The clustering process starts again with a box around the whole mesh. This box is then subdivided \emph{only} along its largest edge, splitting the mesh into two parts. New boxes are drawn around the elements in each part, and the procedure is repeated recursively until the generated clusters have a radius below a certain threshold. In a default setting the same $h_0$ as in Section~\ref{Sec:Box1}  is used for the the threshold at the root level, for the children the division is stopped if the cluster radius is smaller than $h_\ell = 0.5^\ell h_0$. In comparison with the Box1 clustering the number of the subdivisions of the original box is not defined beforehand but determined adaptively to the geometry of the mesh. Numerical experience shows that this approach is computationally cheap, the size of the clusters is better balanced, and very thin clusters can be avoided. 
\subsection{K-Means Clustering}
In this approach the mesh is subdivided iteratively into $n_0$ clusters at the root level using a k-means clustering algorithm. In this study the \texttt{kmeans} function provided by the statistics package of Octave~\cite{octave} was used. Formally the algorithm aims at finding $n_0$ cluster midpoints $\bz_i$, $i = 1, \dots n_0$ that minimize
$$ 
\text{argmin} \sum_{i = 1}^{n_0} \sum_{\bx \in \mathcal{C}_i} ||\bx - \bz_i||^2,
$$
which is usually done by a two step approach. In a first step $n_0$ points $\bz_i$ in space are chosen as potential cluster midpoints and each element midpoint is assigned to its nearest cluster midpoint forming the initial clusters. In a second step the cluster midpoints are updated based on the elements in their respective cluster. This procedure is repeated iteratively until a fixed number of iterations is reached or until the change in cluster midpoints is below a certain threshold. For more information we refer to the \texttt{kmeans} clustering algorithm provided in the statistics package of Octave \cite{octave}.

The k-means clustering is only used for the clustering at the root level, for all other levels of the hierarchical tree we estimate an initial box length $h_0$ and use the child generating process from  either the Box1 or the Box2 clustering.  Motivated by the calculation of the default box length for the Box1 and Box2 clustering, the number of clusters at root level should be about $n_0 \approx 0.9\sqrt{N}$. After the clustering process the initial box lenght estimator $h_0$ is determined as the square root of the mean number of elements in each cluster times the average over the areas of all elements in the mesh divided by 0.9:
$$
h_0 = \frac {\sqrt{ {N}_{\mathcal{C}}  {A_0} }} {0.9}.
$$
\subsection{General Clustering Considerations}\label{Sec:Problems}
One needs to consider that the expansion length of the FMM is defined by the largest radius of the clusters at a level, see, e.g., Eq.~(\ref{Equ:Truncation}). The far-field criterion, on the other hand, is based on the ``local'' radii of the cluster pair involved, see Eq.~(\ref{Equ:FarfieldCrit}). Thus, if there is a mixture of small and large  clusters (in terms of radii) on a level, numerical problems may occur because two small clusters for which the far-field criterion is fulfilled may be still too close for a numerically stable evaluation of the Hankel function of high orders that are necessary due to the high truncation order $L$ ``caused'' by the large clusters at that level.

\begin{figure}[!h]
  \begin{center}
    \begin{tabular}{rcrcrc}
      \raisebox{0.250\textwidth}{\small a) \hspace{-20pt}} &
      \includegraphics[height=0.25\textwidth]{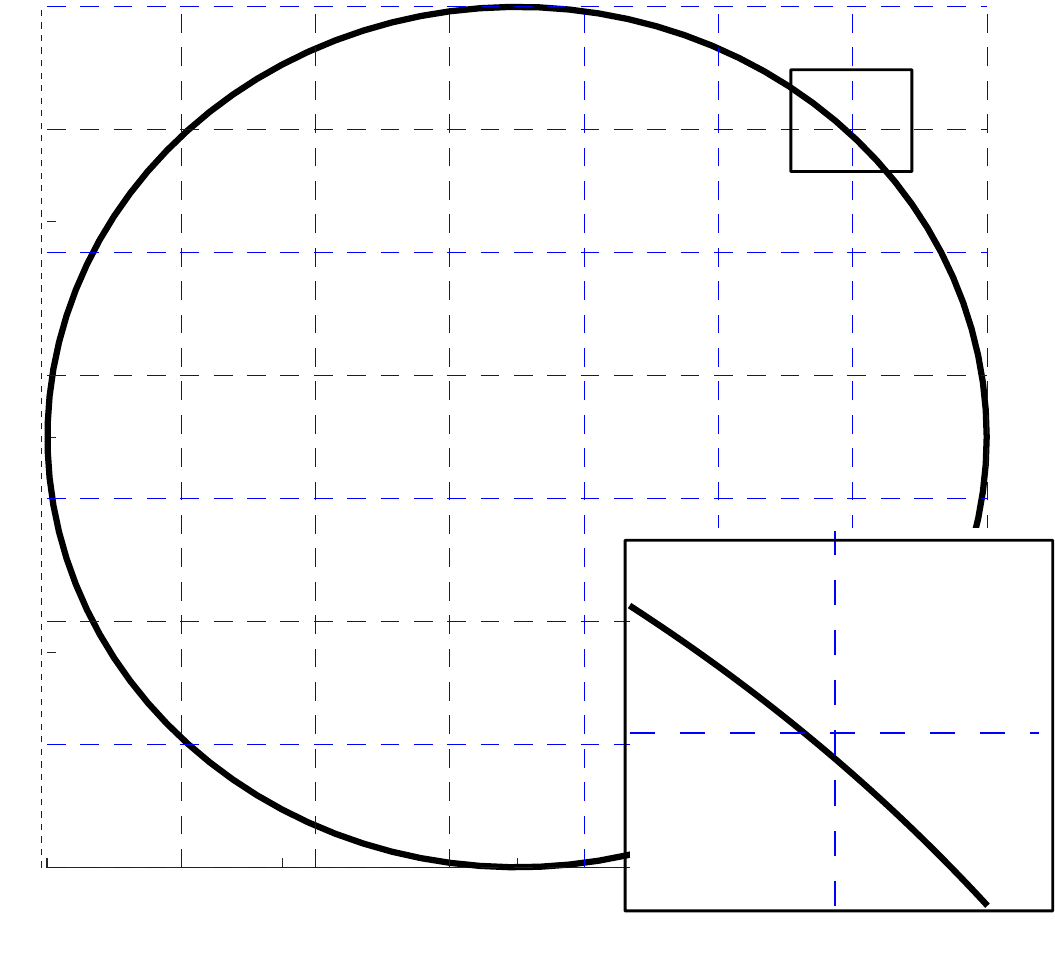}
      \raisebox{0.25\textwidth}{\small b) \hspace{-20pt}} &
      \includegraphics[height=0.2\textwidth]{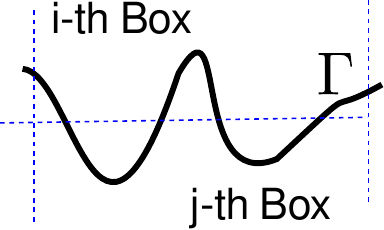}
      \raisebox{0.25\textwidth}{\small c) \hspace{-20pt}} & \includegraphics[height=.25\textwidth]{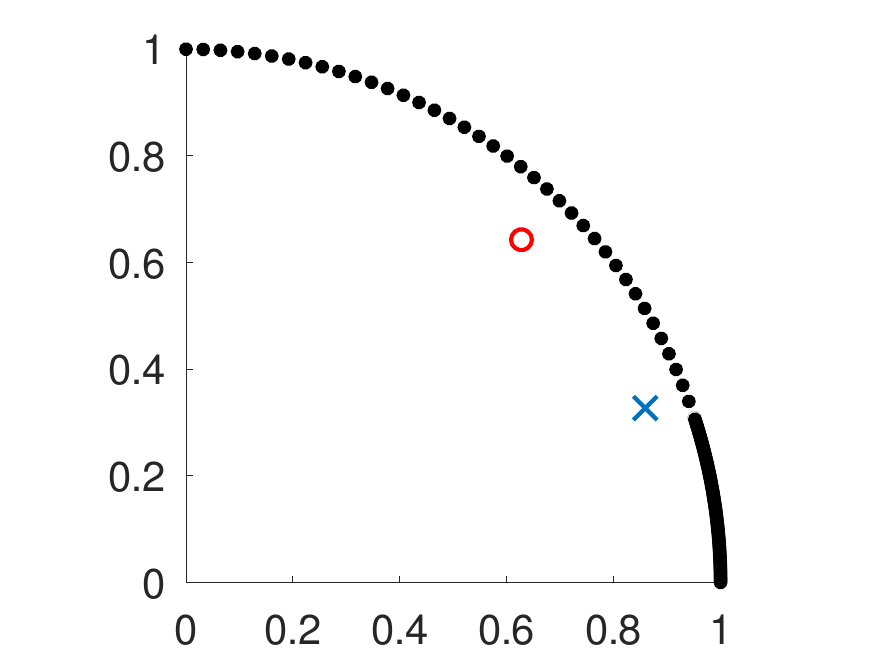}
    \end{tabular}
  \end{center}
  \caption{Example of meshes, where a) some clusters contain only very few elements, and b) elements in the cluster/box are not contiguous. The thick lines depict the surface $\Gamma$ of the object. c) An example, where a mixture of large and small elements inside a single cluster introduce a bias when calculating  the cluster midpoint. 'x' denotes the cluster midpoint calculated as sum over all vertices, 'o' denotes the cluster midpoint that is calculated as a weighted sum over all vertices.}\label{Fig:Problems}
\end{figure}

\begin{figure}[!h]
  \begin{center}
    \includegraphics[width=0.6\textwidth]{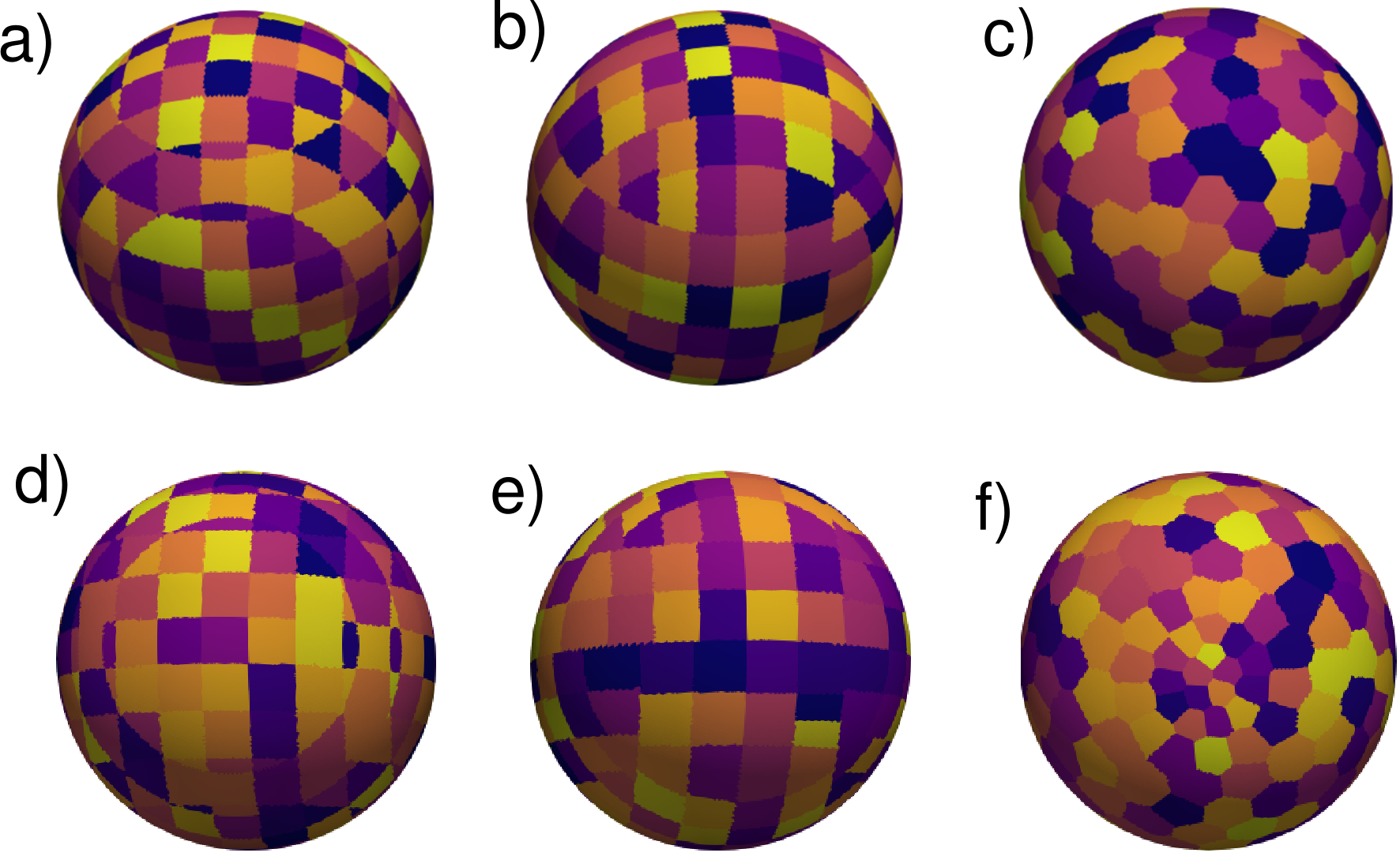}
  \end{center}
  \caption{Example of the three types of clustering for the unit sphere. Top row: a) Box1 clustering, b) Box2 clustering c) k-means clustering for a quasi-uniform mesh. Lower Row: a) Box1 clustering, b) Box2 clustering c) k-means clustering at the north pole. Around the pole the mesh was discretized locally with more elements.}\label{Fig:CompareSphereClustering}
\end{figure}

Even for regular meshes with elements of almost uniform size the Box1 clustering may lead to boxes containing very few elements and/or small radii. One possible scenario is depicted in Fig.~\ref{Fig:Problems}a. The size of the box and the number of subdivisions was chosen such that in the upper right corner a sub-box can be found that contains only a very small portion of the mesh leading to a very small cluster radius.  To avoid this problem, small clusters (again in terms of radii) may be merged with larger neighbors. But even then, it is still possible that in 3D very narrow clusters do occur, see for example Fig.~\ref{Fig:CompareSphereClustering}a. Fig.~\ref{Fig:CompareSphereClustering} shows an example for a clustering at root level using the default $h_0$ settings for a uniform and a non-uniform  mesh of the unit sphere with about 70000 elements. Due to the way of constructing the root of the cluster tree, narrow clusters do not seem to be that prominent in the Box2 clustering.  For the k-means clustering such small clusters can only occur in very extreme cases. 

It is \emph{not} guaranteed that the elements in single cluster form a  contiguous region, see for example Fig.~\ref{Fig:Problems}b. If the geometry is locally not smooth or the mesh includes narrow passages and crests (e.g., the ear canal in a human head model), a box boundary can split the mesh into different parts. In combination with the proposition to merge small clusters with their neighbors it may happen that  children of non-contiguous parents may be merged again into  their parent, which leads to unnecessarily large cluster radii and stability problems. 

If the midpoint of a cluster is calculated as the average over all coordinates of vertices in the cluster, the midpoint will be biased towards regions containing small elements, see Fig.~\ref{Fig:Problems}c. This may lead to unnecessarily large cluster radii. To circumvent this problem a weighting dependent on the area of the elements may be introduced. 

Compared to the box based clustering methods the k-means clustering has the advantage that it is per definition easy to fix the number of clusters beforehand, but it has the  drawback that the clustering depends slightly on the local size of the elements in a region. In Figs.~\ref{Fig:CompareSphereClustering}d) to f) the three clustering methods are compared for a non-uniform mesh of the unit sphere where the area around the north pole (0,0,1) is sampled using finer elements. In the k-means clustering more and smaller clusters appear around this pole whereas in the box-based clustering this effect is not observed.
%%%%%%%%%%%%%%%%%%%%
%
%%%%%%%%%%%%%%%%%%%%
\section{Numerical Experiments}\label{Sec:Experiments}
For the numerical experiments in this section a modified version of the open source code \numcalc{} is used \cite{Kreuzeretal24}. Experiments are performed using two different geometries. First, we look at the classical benchmark problem of the acoustic field on the sound hard unit sphere caused by a plane wave. For this problem we will use 3 different discretizations, one where the elements have about the same size (see Section~\ref{Sec:Uni} for a description), and two non-uniform meshes where the discretization around the north pole of the sphere gets finer, see Section~\ref{Sec:NonUni} for a description of these non-uniform meshes. For these three meshes we will look at the structure of the cluster tree for different initial box lengths and the three different clustering methods. We will compare the computation time and the memory consumption as a function of initial box lenght. Note that a direct comparison of the different clustering methods just based on $h_0$ is not trivial because $h_0$ has a slightly different meaning for the different clustering methods. For the Box1 clustering $h_0$ is actually the edge length of a box, but it is modified because of rounding, see Section \ref{Sec:Box1}. For the Box2 clustering $h_0$ is a threshold for the radius of a cluster. If the radius is above $h_0$, the cluster will be split into two parts. For the k-means clustering $h_0$ is determined \emph{after} the clustering in the root level and acts as an auxiliary variable to determine the clusters for the rest of the tree.

The second type of geometry investigated will be real life problem of calculating head related transfer functions (HRTFs) using the boundary element method. HRTFs are functions that model the filtering effect of the human head and the human pinna on incoming sound. This effect is used by humans to determine sound source positions in 3D \cite{Moeller92}. To that end we will look at two meshes of a human head, where the outer ear (= pinna) is discretized much finer then the rest of the head. We will directly compare the performance of the three different meshing methods proposed in Section \ref{Sec:Clusteringmethods}.
\subsection{Sphere Meshes}
\subsubsection{Uniformly Meshed Sphere}\label{Sec:Uni}
The uniform discretization of the unit sphere  is based on an 20-sided icosahedron, where each side is subdivided into $n^2$ similar subtriangles where $n$ is a number provided by the user. The vertices of these subtriangles are then projected onto the unit sphere.

\begin{figure}[!h]
  \includegraphics[width=0.48\textwidth]{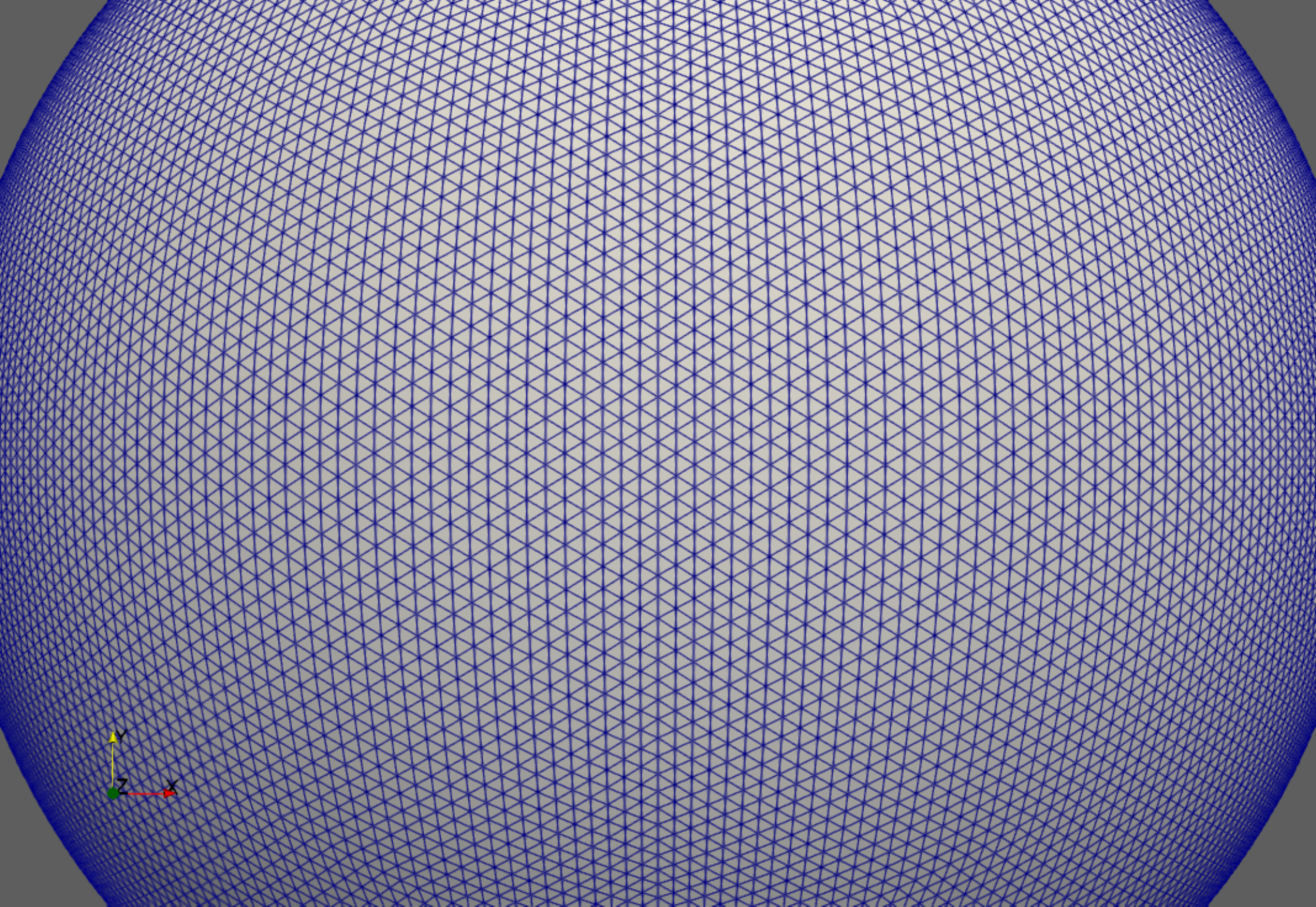}
    \includegraphics[width=0.48\textwidth]{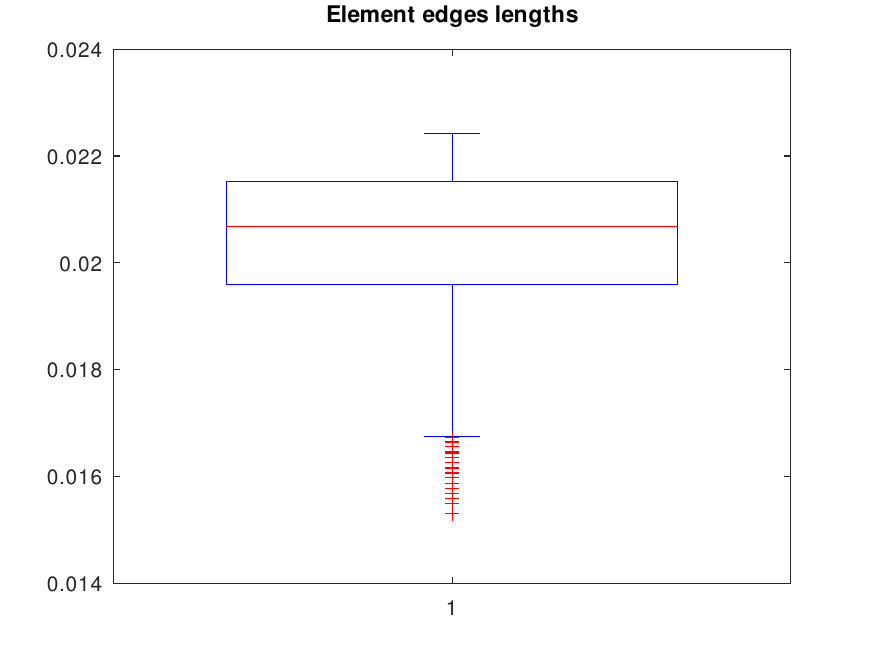}
  \caption{Uniform meshing and boxplot of the edge lengths for  the mesh modeling the unit sphere.}\label{Fig:RegSphere}
\end{figure}

The uniform mesh consists 69620 triangular elements and 34812 vertices. The number of elements was chosen such that the numerical errors of the BEM up to 2000 Hz are in an acceptable range. In Fig.~\ref{Fig:RegSphere} parts of the mesh and a boxplot of the edge lengths of the elements contained in the mesh are depicted. Although there are a few small elements compared to the average element, the mesh can be considered as regular. 

\subsubsection{Non-Uniform Sphere-Meshes}\label{Sec:NonUni}
\begin{figure}[!h]
  \begin{tabular}{rcrcrc}
    \raisebox{0.250\textwidth}{\small a) \hspace{-20pt}} &
    \raisebox{8pt}{\includegraphics[width=0.29\textwidth]{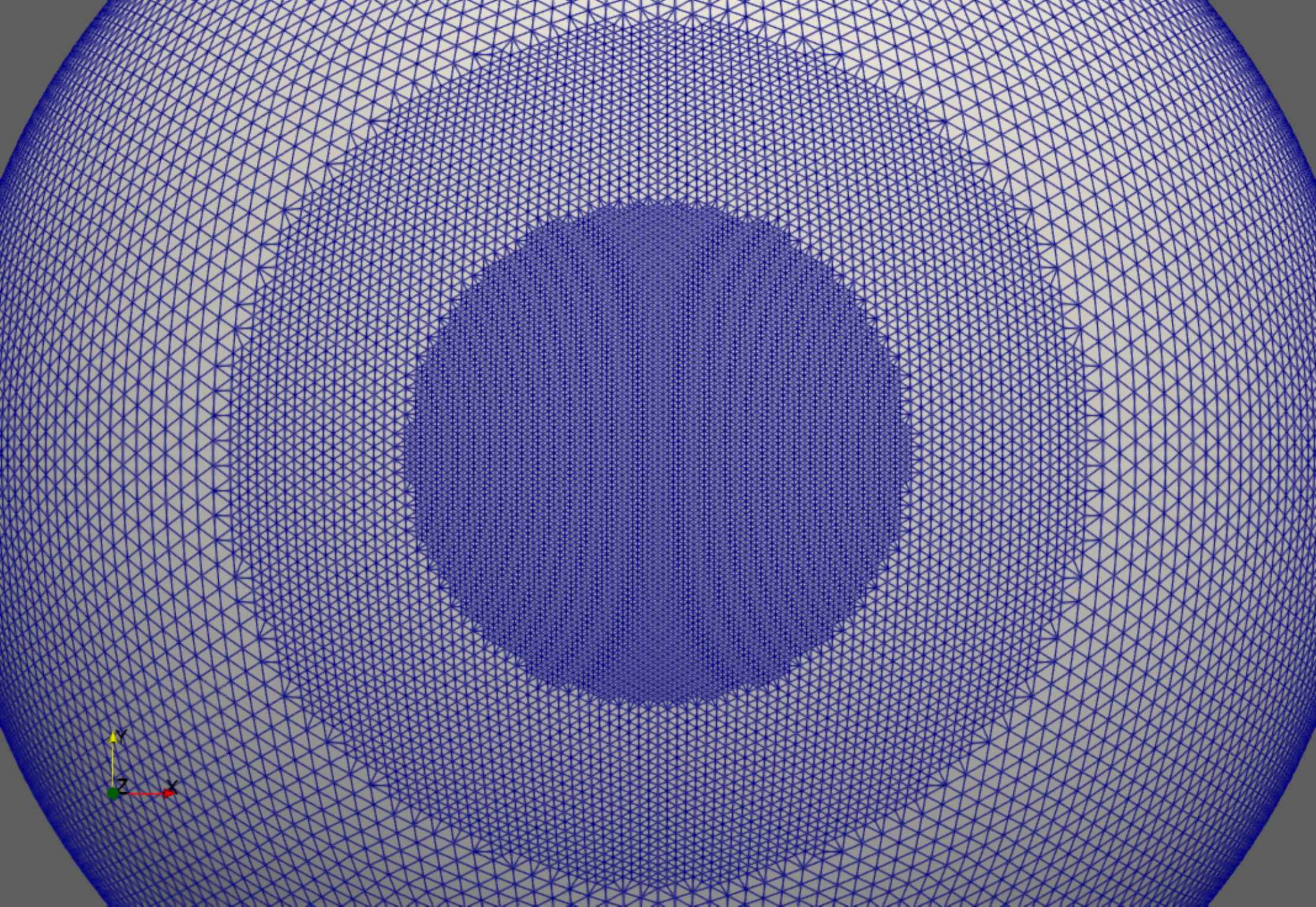}} &
    \raisebox{0.250\textwidth}{\small b) \hspace{-20pt}} &
    \raisebox{8pt}{\includegraphics[width=0.29\textwidth]{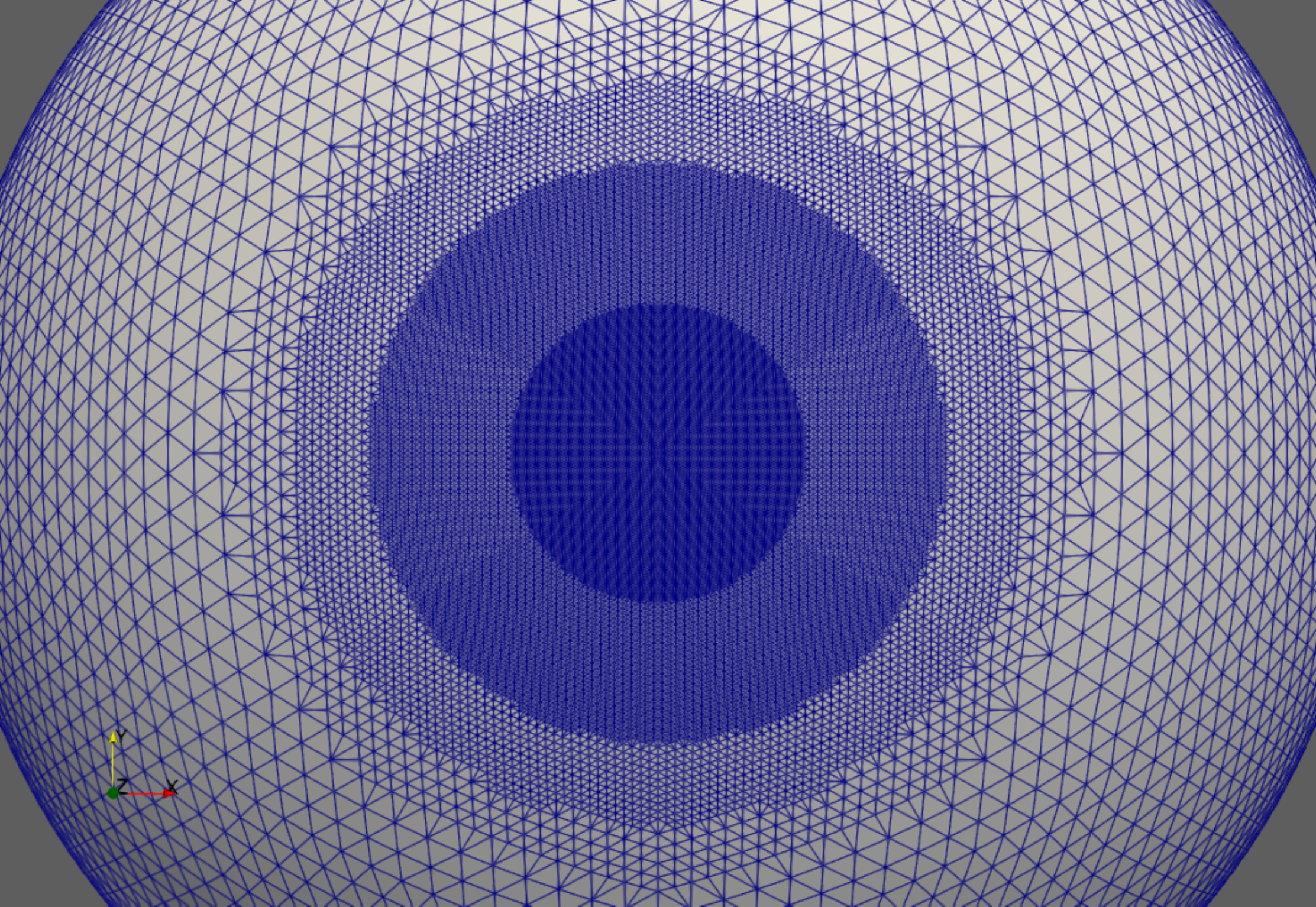}} &
    \raisebox{0.250\textwidth}{\small c) \hspace{-20pt}} &
    \includegraphics[width=0.3\textwidth]{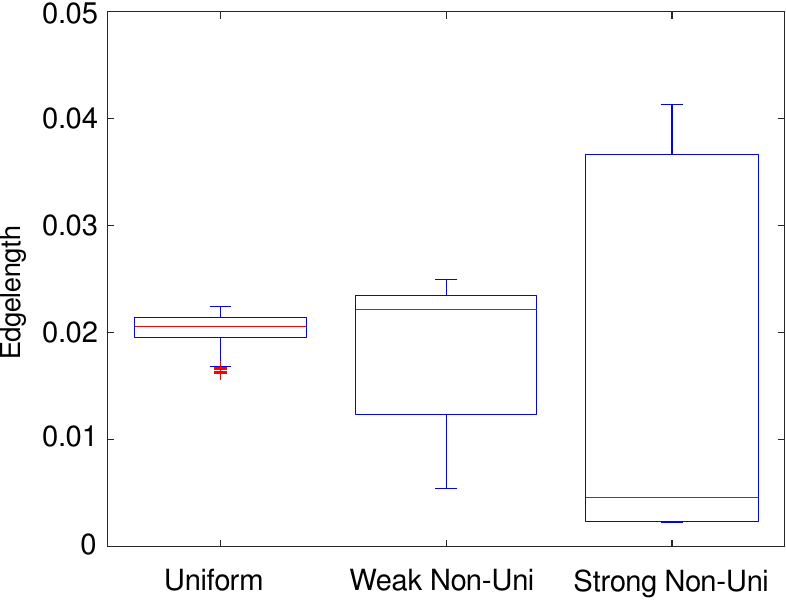}
  \end{tabular}
  \caption{a) Weak Non-Uni  and b) Strong Non-Uni meshes  around the north pole and c) boxplots of the edge lengths of each triangle in the spherical meshes meshes.}\label{Fig:NonUniformSpheres}
\end{figure}

To evaluate possible problems caused by non-uniform meshes we look at spherical meshes where the element sizes around the north pole is smaller compared to the rest of the sphere. In the first mesh (denoted as Weak Non-Uni mesh in the following) about 29\% of all vertices in the mesh have a $z$-coordinate bigger than 0.85. The second non-uniform mesh (Strong Non-Uni) is slightly coarser farther away from the north pole but the elements around the north pole have been further refined leading to a mesh with 73\% of the vertices have a $z$-coordinate bigger than 0.8. Fig.~\ref{Fig:NonUniformSpheres} illustrates the discretization around the north pole for the Weak Non-Uni mesh (left), the Strong Non-Uni mesh (center) and boxplots comparing the distribution of edge lengths of each sphere-mesh (right).

\subsection{Clustering Information Sphere-Meshes}
In this section we look at the dependence of the clustering as a function of the initial box length $h_0$ used to create the root cluster with the three different clustering methods. For the spherical meshes calculations were done up to the initial box length of $h_0 = 0.85$\,m. At this value the root of the cluster tree for the uniform mesh and the Box1 clustering contains only 8 clusters. No multipole expansion is possible at this level, because no cluster pair fulfills the far field criterion given in Eq.~(\ref{Equ:FarfieldCrit}).
\begin{figure}[!h]
  \begin{tabular}{r}
    \multicolumn{1}{c}{Uniform}\\
    \includegraphics[width=0.96\textwidth]{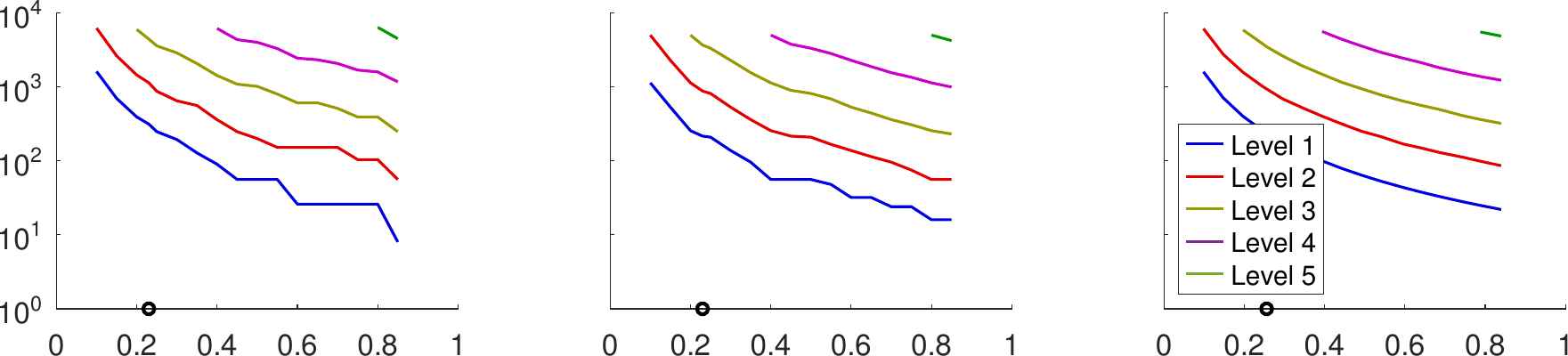}\\
    %    Non-Uni1
    \multicolumn{1}{c}{Weak Non-Uni}\\
    \includegraphics[width=0.98\textwidth]{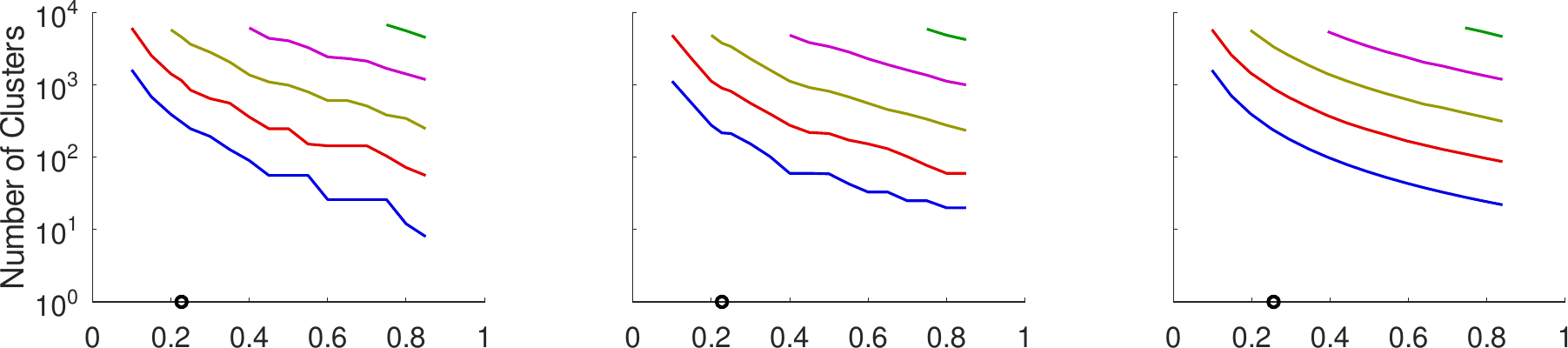}\\
    %    Non-Uni2
    \multicolumn{1}{c}{Strong Non-Uni}\\
    \includegraphics[width=0.96\textwidth]{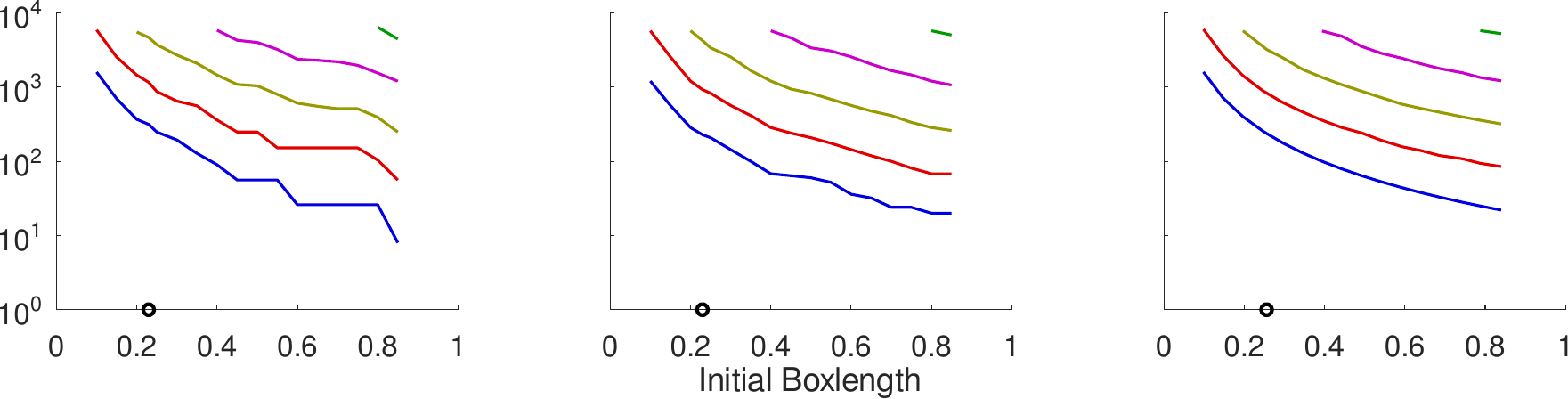}
  \end{tabular}
  \caption{Number of clusters in each level as a function of the initial box length $h_0$. The first row depicts the data for the Uniform mesh, the second row for the Weak Non-Uni mesh, and the third row the data for the Strong Non-Uni mesh. The first column depicts the data for the Box1 clustering, the second column for the Box2 clustering and the third column for the k-means clustering. The circles represent the default $h_0$ for the three clustering methods, see Eq.~(\ref{Eq:DefaultBox})}\label{Fig:CompNclustall}
\end{figure}
\begin{figure}[!h]
  \begin{tabular}{r}
    \multicolumn{1}{c}{Uniform}\\
    \includegraphics[width=0.96\textwidth]{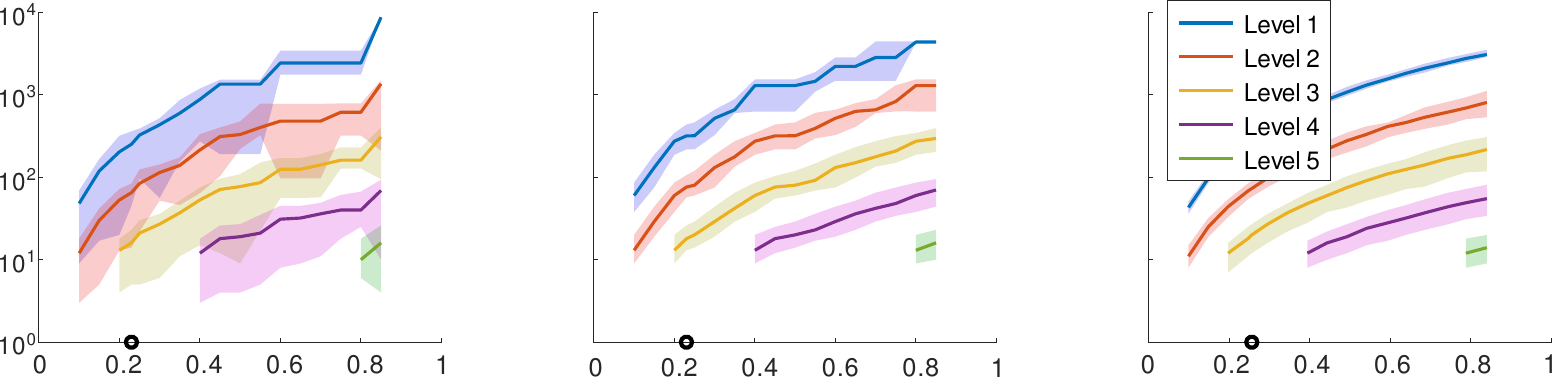}\\
    \multicolumn{1}{c}{Weak Non-Uni}\\
    \includegraphics[width=0.98\textwidth]{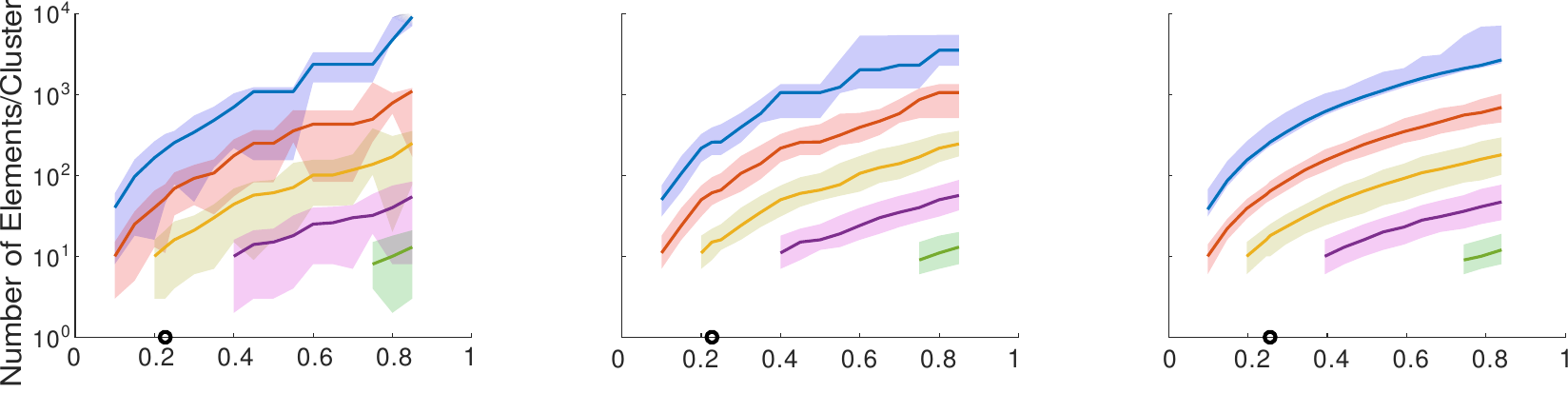}\\
    \multicolumn{1}{c}{Strong Non-Uni}\\
    \includegraphics[width=0.98\textwidth]{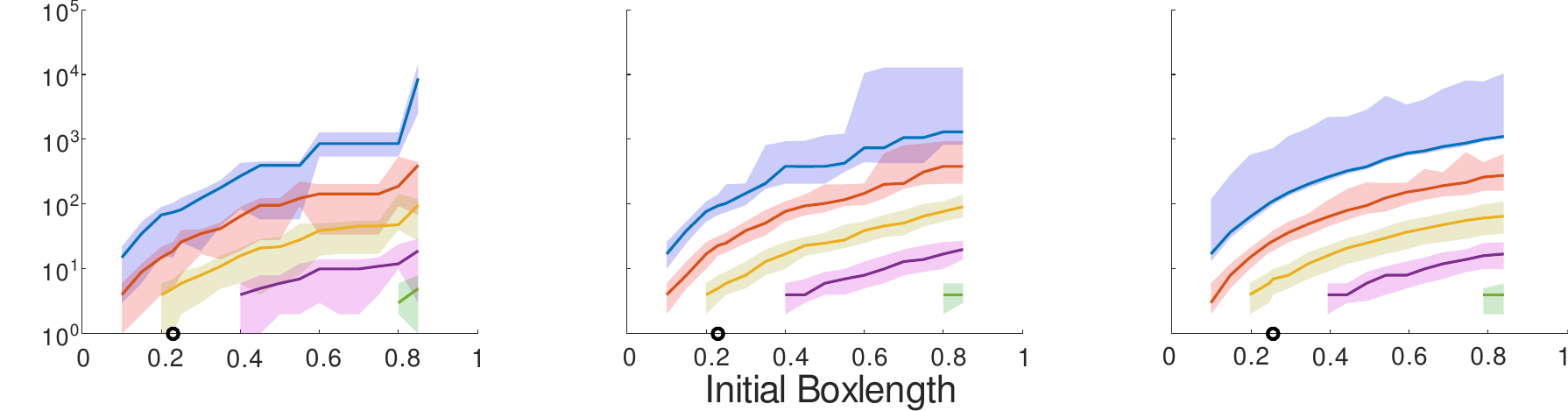}
  \end{tabular}
  \caption{Median number of elements in the clusters for each level as a function of the initial box length. The first row depicts the data for the Uniform mesh, the second row for the Weak Non-Uni mesh, and the third row the data for the Strong Non-Uni mesh. The first column depicts the data for the Box1 clustering, the second column for the Box2 clustering and the third column for the k-means clustering. The circles represent the default $h_0$ for the three clustering methods. The shaded areas depict the range between  lower and upper 10\% quantile. }\label{Fig:MedianSphere}
\end{figure}

In Fig.~\ref{Fig:CompNclustall} the number of clusters in each level for all three meshes (rows: Uniform, Weak Non-Uni, and Strong Non-Uni) and the three clustering methods (columns: Box1, Box2, and k-means clustering) is depicted Fig.~\ref{Fig:MedianSphere} shows the median value and the 10\% percentiles of the number of elements in each cluster as a function of initial box length for all three meshes and all three  clustering methods. 
Already for the uniform mesh a relatively big variance in the number of elements per cluster can be observed for the Box1 clustering. This can be explained by the fact that for the given initial box length some boxes still may can contain very few elements, see Fig.~\ref{Fig:Problems}a for an example. Merging small (in terms of cluster radius) clusters with a neighboring cluster does not prevent clusters to become very slim. As the cluster radius is defined as maximum distance from cluster midpoint to its vertices, slim clusters will have a ``normal'' radius and, therefor, not merging will be done.

\begin{figure}[!h]
  \begin{tabular}{r}
    \multicolumn{1}{c}{Uniform}\\
    \includegraphics[width=0.98\textwidth]{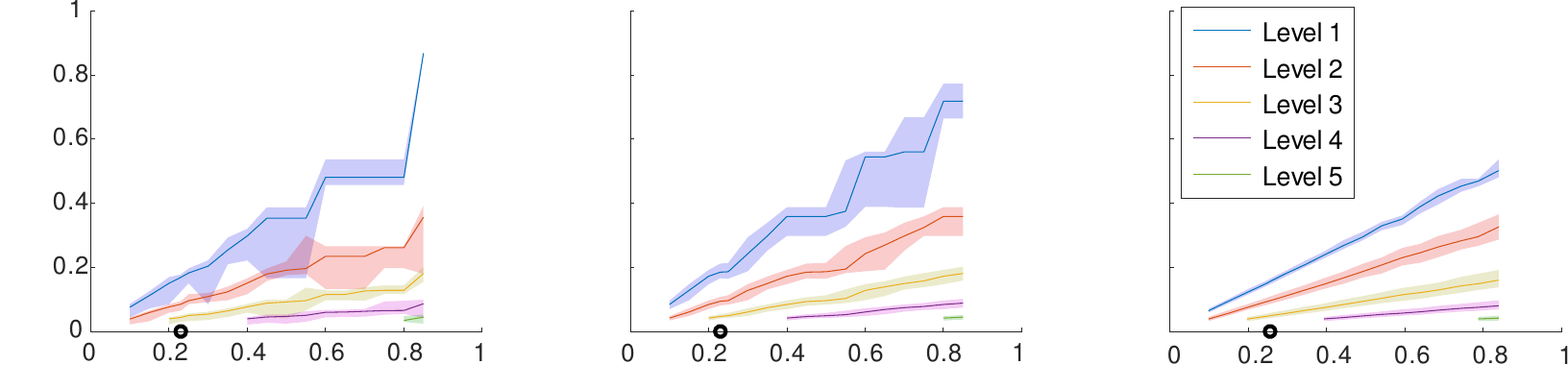}\\
    \multicolumn{1}{c}{Weak Non-Uni}\\
    \includegraphics[width=0.98\textwidth]{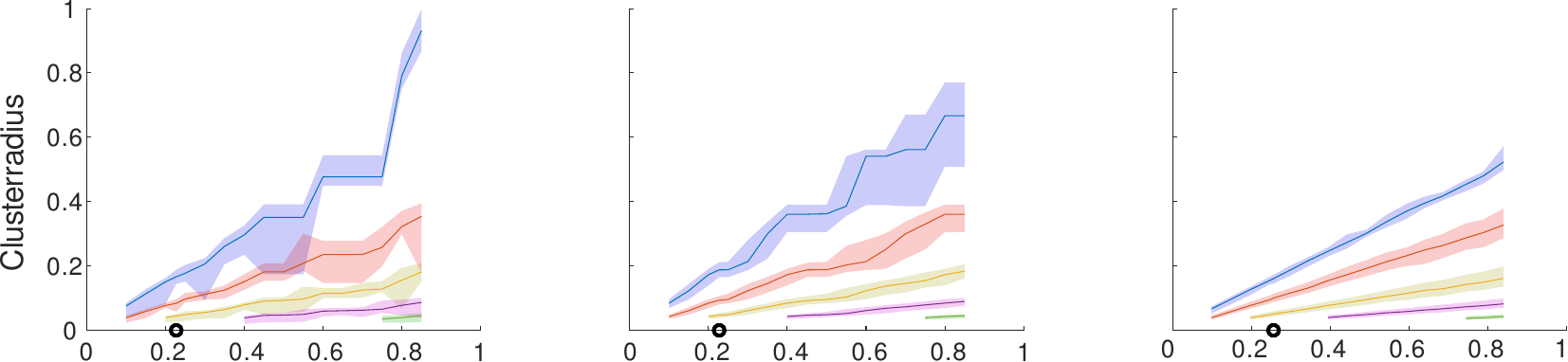}\\
    \multicolumn{1}{c}{Strong Non-Uni}\\
    \includegraphics[width=0.98\textwidth]{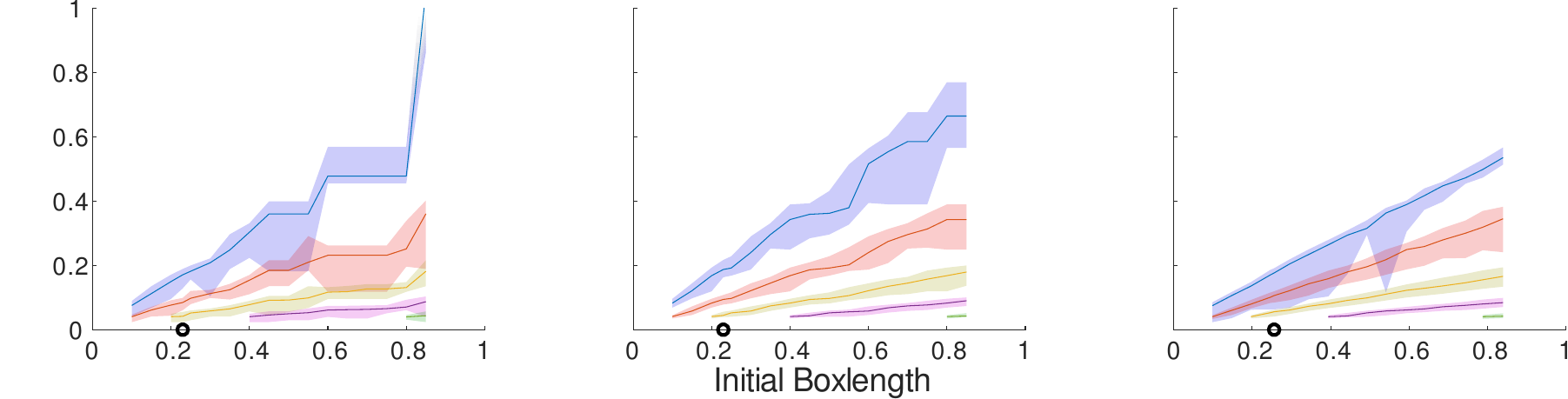}
  \end{tabular}
  \caption{Median number for the radii of the clusters for each level as a function of the initial box length. The first row depicts the data for the Uniform mesh, the second tow for the Weak Non-Uni mesh, and the third row the data for the Strong Non-Uni mesh. The first column depicts the data for the Box1 clustering, the second column for the Box2 clustering and the third column for the k-means clustering. The circles represent the default $h_0$ for the three clustering methods. The shaded areas depict the range between  lower and upper 10\% quantile.}\label{Fig:ComprSphere}
\end{figure}

This effect is reduced for the Box2 clustering, and even more for the k-means clustering. When looking at the different radii of the clusters in Fig.~\ref{Fig:ComprSphere} a similar effect can be observed. In the figures it can also be observed that for the Box1 clustering, the number of elements per cluster at the root level is the same for multiple initial box lengths $h_0$. This is because the the actual edge length of the sub boxes in the Box1 clustering is determined using the \emph{rounded} value, see Section \ref{Sec:Box1}. 

For the non-uniform meshes the results indicate that the variation in number of elements per cluster and the cluster radii gets larger for the k-means clustering, which can be explained by the fact that the sizes of the k-means clusters vary locally with the size of the elements. This effect is reduced in higher levels because a parent is only split into children, if its radius is big enough.

For the non-uniform meshes of the unit sphere there are extreme cases where the cluster radii can become bigger than 1. This is because the cluster midpoint is calculated as the mean value over all vertices in the cluster, thus, many small clusters introduce a bias towards their vertices, see for example Fig.~\ref{Fig:Problems}c.
\subsubsection{Expansion length}
\begin{figure}[!h]
  \begin{center}
  \includegraphics[width=0.48\textwidth]{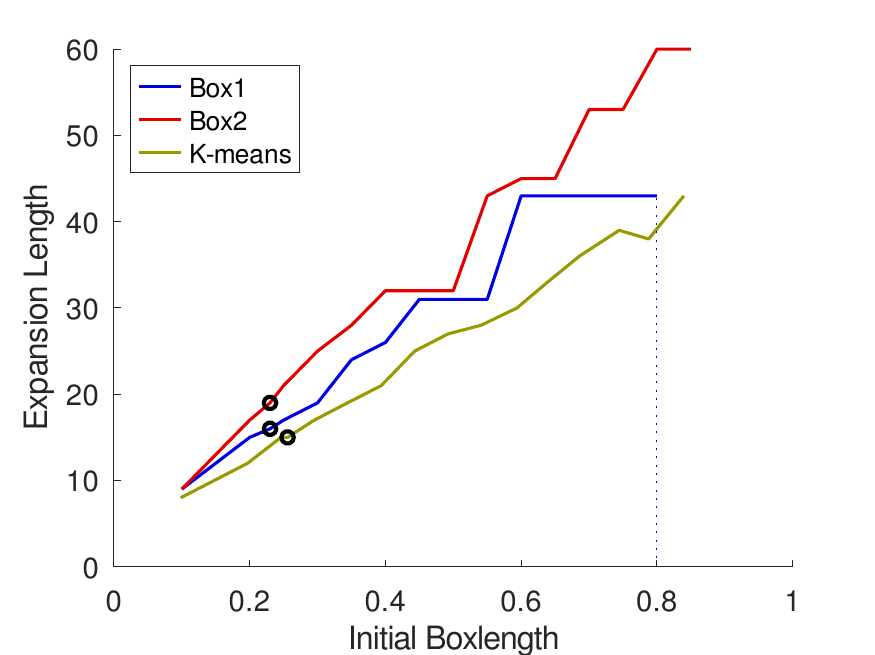}
  \includegraphics[width=0.48\textwidth]{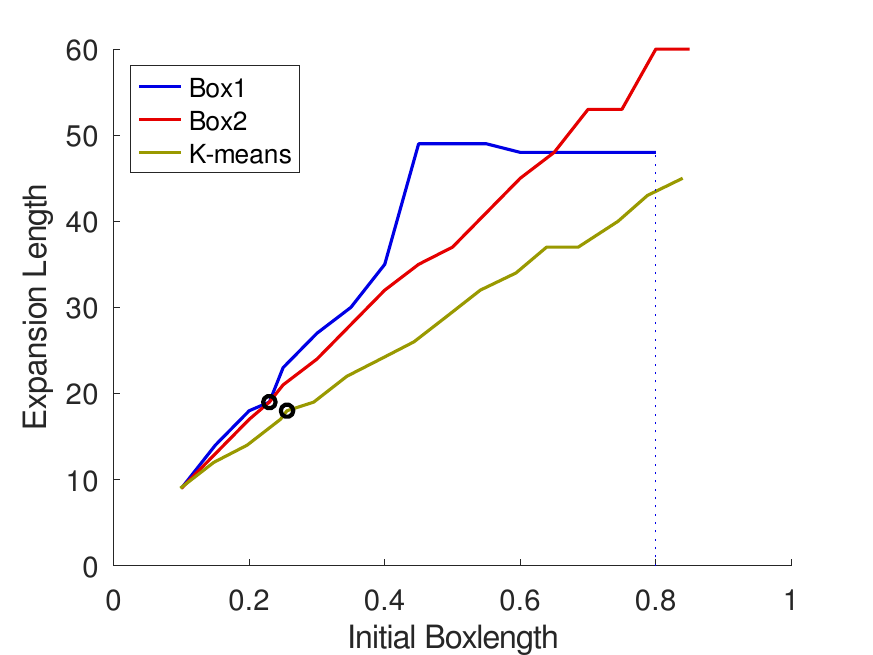}
  \end{center}
  \caption{Expansion length for the different clustering method at 2000~Hz at the root level (Level 1) as a function of initial box length. The graph on the left shows the expansion length for the uniform mesh, the right graph the expansion length for the Strong Non-Uni mesh. The circles indicate the default box length for each clustering version.}\label{Fig:Explengths0}
\end{figure}

\begin{figure}[!h]
  \begin{center}
  \includegraphics[width=0.48\textwidth]{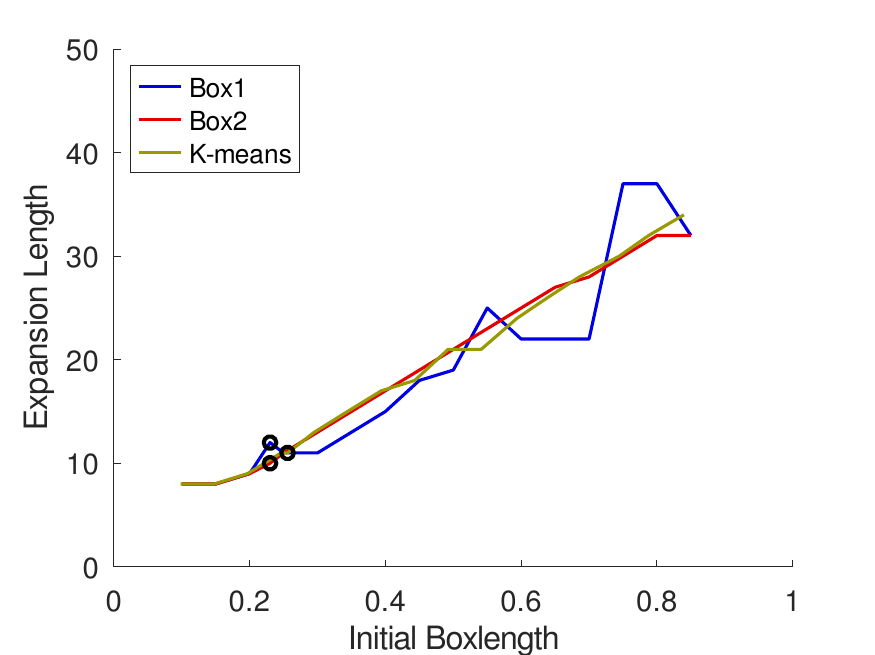}
  \includegraphics[width=0.48\textwidth]{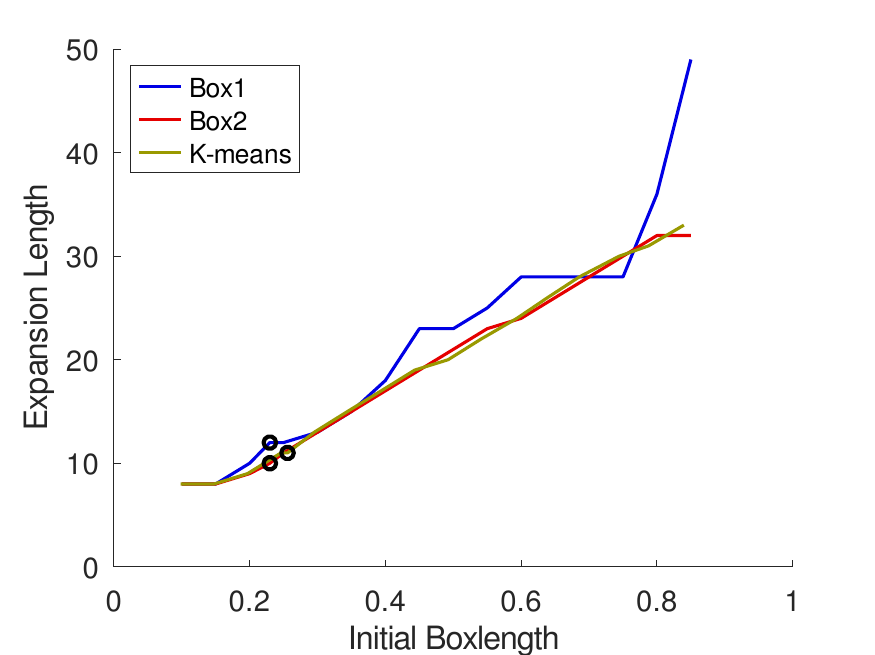}
  \end{center}
    \caption{Expansion length for the different clustering method at 2000~Hz at Level 2 of the MLFMM  as a function of initial box length. The graph on the left shows the expansion length for the uniform mesh, the right graph the expansion length for the Strong Non-Uni mesh. The circles indicate the default box length for each clustering version.}\label{Fig:Explengths1}
\end{figure}

In Figs. \ref{Fig:Explengths0} and \ref{Fig:Explengths1} the expansion lengths in the root (first multipole level) and second multipole level are depicted at a frequency of 2000~Hz for the uniform and the Strong Non-Uni mesh. For the non-uniform mesh 2000~Hz is the upper frequency limit for which the BEM calculations still yields  results with acceptable accuracy. For lower frequencies the curves look similar except that they are shifted downwards. For $h_0 = 0.85$\,m the Box1 clustering has only 8 clusters at the root level, and there are no far field cluster pairs. Thus, the expansion length is set to 0. At the root level there is some influence of the clustering method on the the expansion lengths, but these differences almost vanishe at higher multipole levels. Numerical experience has shown that expansion lengths above 50 may become problematic because standard implementations of the Hankel function can exhibit stability problems for such high orders. Nevertheless, for the default box lengths the respective expansion lengths are well within a ``safe'' range. 
\subsubsection{Performance}
\begin{figure}[!h]
  \includegraphics[width=0.98\textwidth]{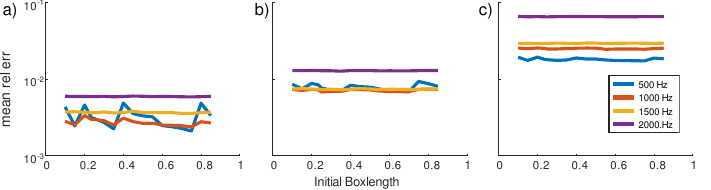}
  \caption{Relative mean error for a) the uniform, b) the Weak Non-Uni, and c) the Strong Non-Uni spherical mesh at frequencies 500, 1000, 1500, and 2000 Hz using the MLFMM with the Box1 clustering.}\label{Fig:Err}
\end{figure}

In Fig.~\ref{Fig:Err} the mean relative error of the sound pressure over all element midpoints for the three different spherical meshes is depicted as a function of the initial box length $h_0$. As there is almost no difference in the error between the three different clustering approaches, only the Box1 clustering results are depicted. It is apparent that the \emph{overall} error for the non-uniform meshes still depends on the coarse discretization away from the north pole of the sphere. It can be observed that locally the error at the finer discretized parts of the meshes is slightly smaller, however, for the elements in the border region between small and large elements the errors become relatively large. This behavior is per se not caused by MLFMM as it can also be observed for the  traditional BEM approach. Note that this error  does \emph{not} contradict the observations made in \cite{Ziegelwangeretal16} dealing with non-uniform meshes for HRTF calculations. First, it was pointed out in \cite{Ziegelwangeretal16} that the error can be \emph{perceptually} neglected, second, the meshes used in \cite{Ziegelwangeretal16}  were refined \emph{continuously} towards the region of interest, in this case the pinna and the ear canal. Third, for the spherical meshes used in this study the mesh is refined at the north pole, that lies in the shadow region of the plane wave, where the absolute sound pressure is very small to begin with. 

\begin{figure}[!h]
  \includegraphics[width=0.5\textwidth]{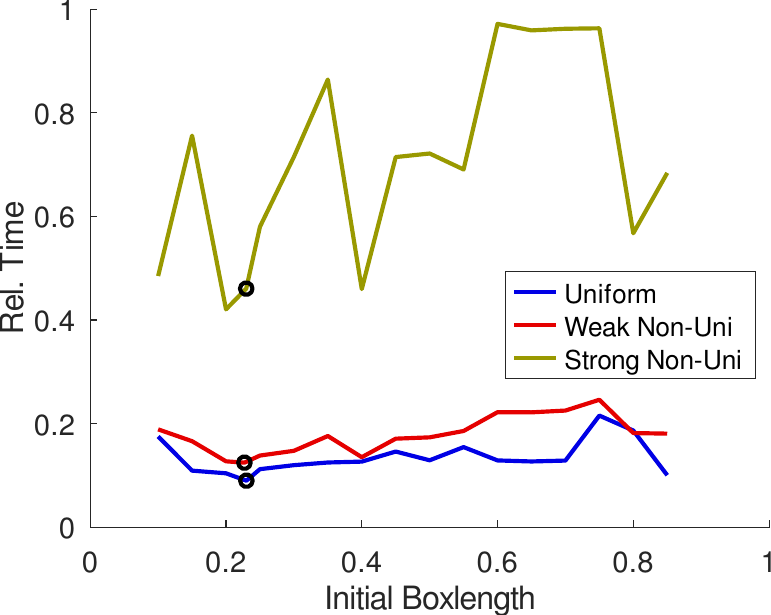}
  \includegraphics[width=0.5\textwidth]{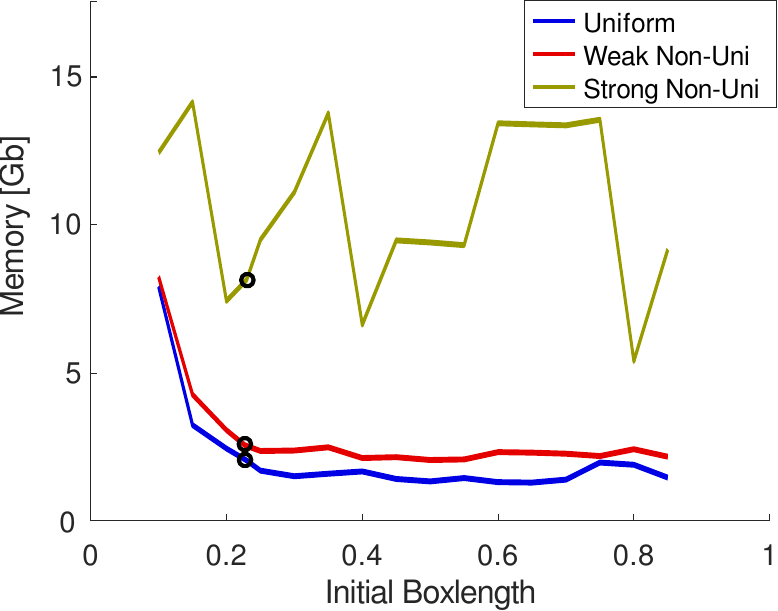}
  \caption{Relative computation time and memory requirements in Gb for the uniform mesh ($N = 69620$ triangular elements), Weak Non-Uni mesh ($N = 73226$), and Strong Non-Uni mesh ($N = 69896$). The black dots denote the position of the default initial box length. In each run the field for 4 frequencies were calculated}\label{Fig:WtimeMem}
\end{figure}

When looking at the computation times it becomes clear, that the non-uniform meshing introduces some additional computational effort. The computation times depicted in  Fig.~\ref{Fig:WtimeMem} are the times for a calculation of the solution for 4 frequencies ($f = 500, 1000, 1500$, and $2000$ Hz) based on one \emph{single} CPU-thread of the same PC to have a better comparison between the different implementations. As the computation time depends on the CPU used, all values are in relation to the maximum time needed for the calculations, in this case 3500\,s. Notches in Fig.~\ref{Fig:WtimeMem} for the Strong Non-Uni mesh correspond to box sizes where an additional level was added to the multi-level approach ($h_0 \approx 0.2,0.4,0.8$\,m, see  Fig.~\ref{Fig:MedianSphere}).

The computation time can be mainly explained by the number of elements in the leaf clusters, which in turn influences the non-zero entries in the nearfield matrix. If we look for example at the results for the Strong Non-Uni mesh at $h_0 = 0.35$\,m about 92\% of the 2092 clusters in the leaf level contain less then 17 elements, however there are some outliers. 12 Cluster at leaf level contain more then 1000 elements. This implies that although these clusters have very small radii, they heavily increase the number of non-zeros in the nearfield matrix, and thus the computation time and memory consumption.

Fig~\ref{Fig:WtimeMem} also depicts the peak memory consumption as a function of the initial box length. Note that contrast to the description given in  \cite{Kreuzeretal24} a modified version of \texttt{NumCalc} was used where the local expansions are  interpolated and filtered between levels. 
\subsection{HRTF Meshes}\label{Sec:HRTFs}
As a practical example of an irregular mesh, a mesh of the human head was chosen that is used to calculate head related transfer functions (HRTFs).
Over the years a multitude of meshes have been created in connection with the Mesh2HRTF project \cite{Brinkmannetal23,Kreuzeretal24,Ziegelwangeretal15} or by other projects~\cite{WangYu26} either by scanning a human head or by generating virtual pinnae and heads. A lot of these meshes were constructed under the assumption that while the pinna needs to be modeled with high accuracy, the details for the rest of the head can be neglected, leading to non-uniform meshes that are fine around the pinnae and coarse at the rest \cite{Ziegelwangeretal16}. One reason for such small elements at the pinna is the fact that many applications use HRTFs up to a frequency of 20\,kHz. Thus, in order to follow the six-to-eight elements per wavelength rule, elements with edge lengths of about 2\;mm are common. A second reason for using such small elements for the pinna is its non-smooth geometry. To avoid numerical problems caused by thin structures, the element sizes are restricted. On the other hand, there is the conception that fewer elements imply less computational effort and memory requirements, thus, the use of bigger elements for the rest of the head. 

\begin{figure}[!h]
  \includegraphics[height=0.3\textwidth]{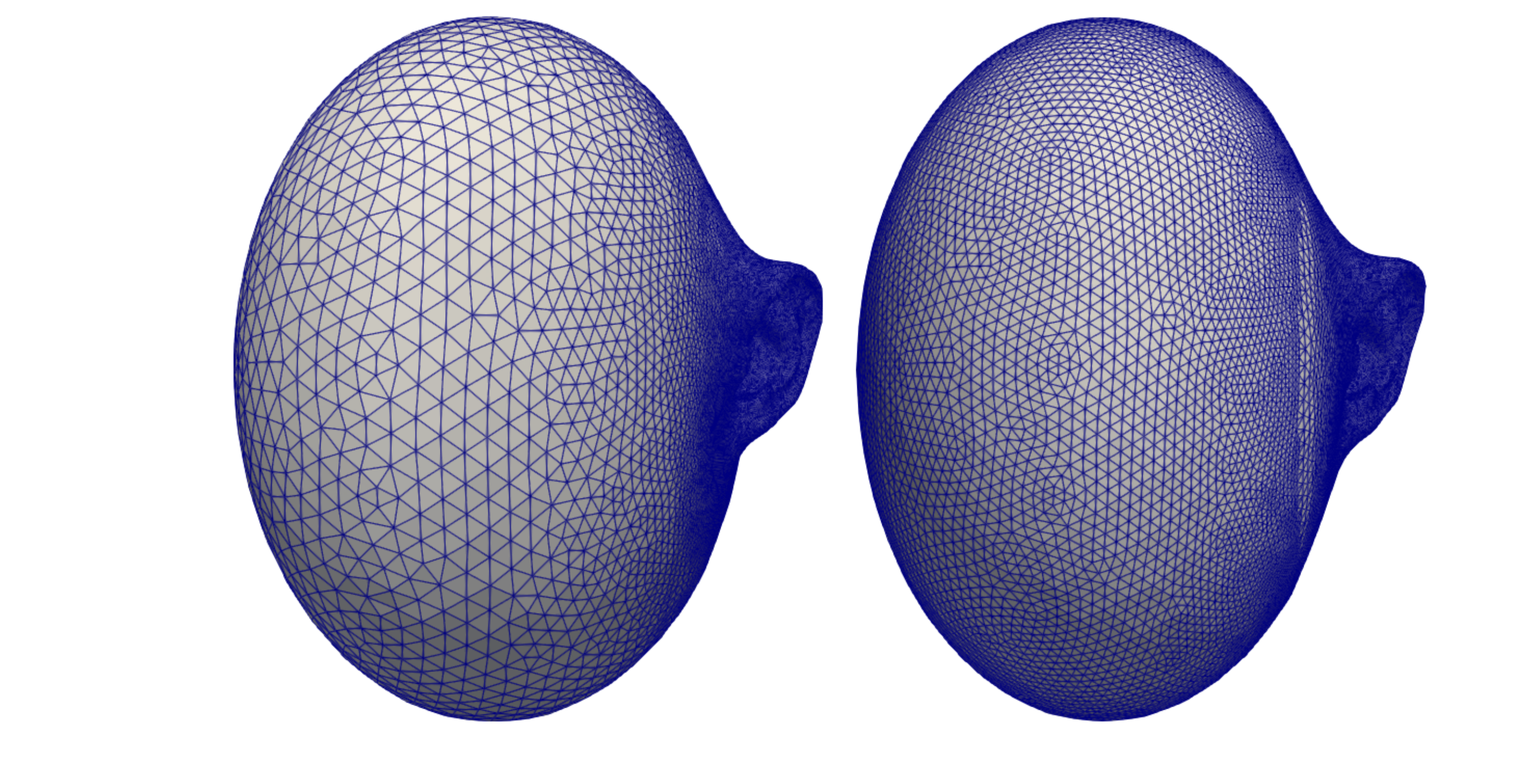}
  \includegraphics[height=0.3\textwidth]{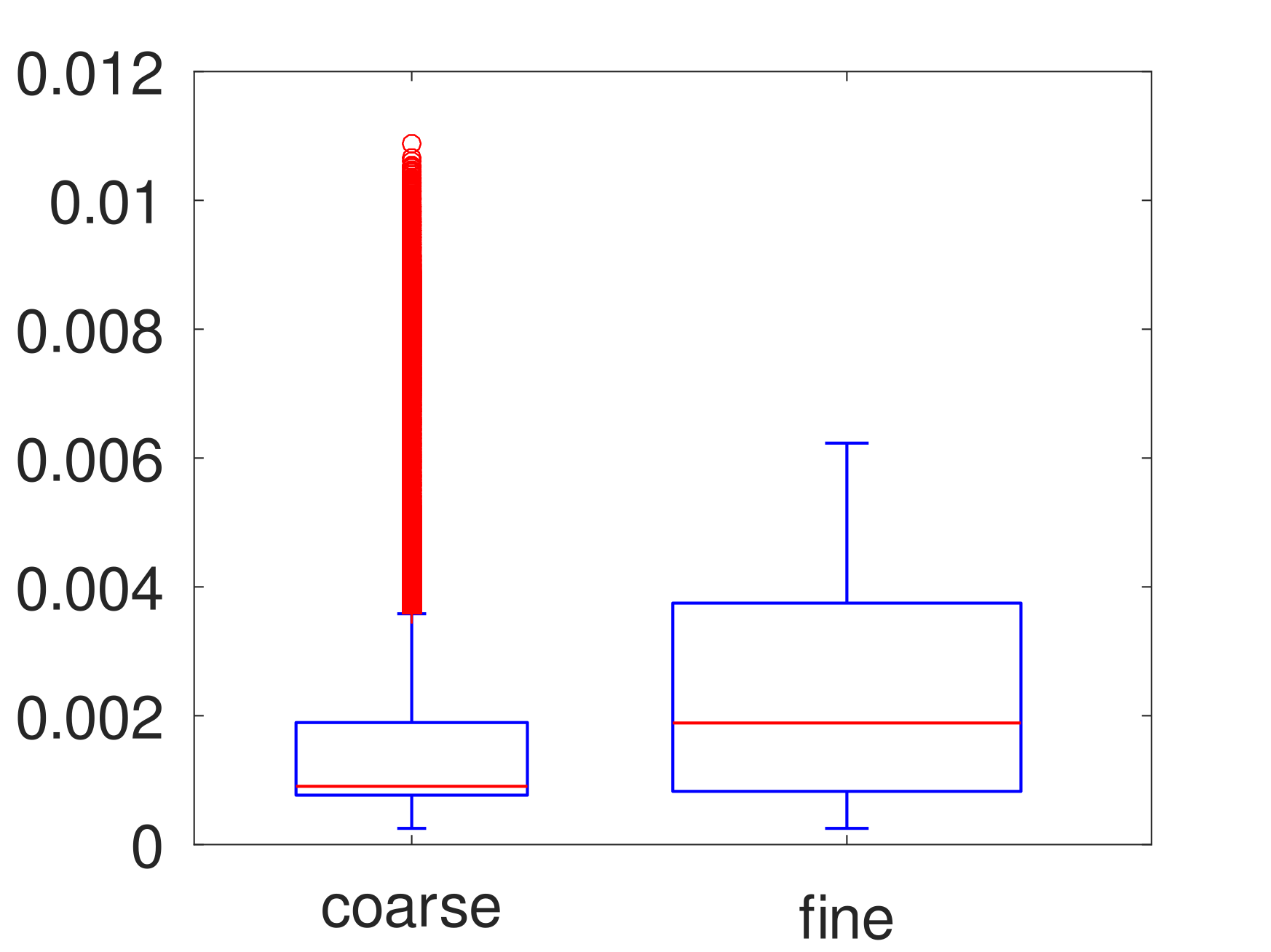}
  \caption{Coarse and fine head mesh and a standard boxplot containing median and quartiles with respect to the edge lengths for the coarse (N = 31050)  and fine mesh (N = 47718)}\label{Fig:Lukasmesh}
\end{figure}

Most of these meshes can be easily treated by \numcalc, but with some meshes like the coarse head mesh  presented in this section numerical problems were encounter, which will be investigated and explained in detail in the following. An additional aim of this section is the comparison of the coarse ``original'' mesh and a finer version  where large elements were split once. The coarse original mesh contains $N = 31050$ elements and was used to calculate HRTFs for frequency up to 22 kHz. For this mesh about 80\% of its elements lie in the area around the pinna. The second, finer mesh was created by subdividing large elements of the original head model and contains $N = 46718$ elements. Both meshes and the boxplots for their respective edge lengths are depicted in Fig.~\ref{Fig:Lukasmesh}. For the calculations a modified version of \numcalc{} was used, where the data between levels was interpolated/filtered. 
%  \begin{figure}[!h]
%    \includegraphics[width=0.48\textwidth]{Pics/Lukasmeshcoarse.pdf}
%    \includegraphics[width=0.48\textwidth]{Pics/Lukasmeshedges.pdf}
%    \caption{Example of a head mesh and the boxplot with respect to the edge lengths of the elements used in the mesh.}\label{Fig:Head}
%  \end{figure}

% mean = 2.12e-3, var = 6.13e-6  var/max = 5.6e.4, var/mean = 2.8e-3
% mean = 2.23e-3, var = 2.12e-6, var/max = 3.4e-4, var/mean = 9.5e-4
%
%
%
%%%%%%%%%%%%%%%%%%%%%%%%%%%%%%%%%%%%%%%%%%%%%%%
\subsubsection{Cluster information}
\begin{figure}[!h]
  \includegraphics[width=0.98\textwidth]{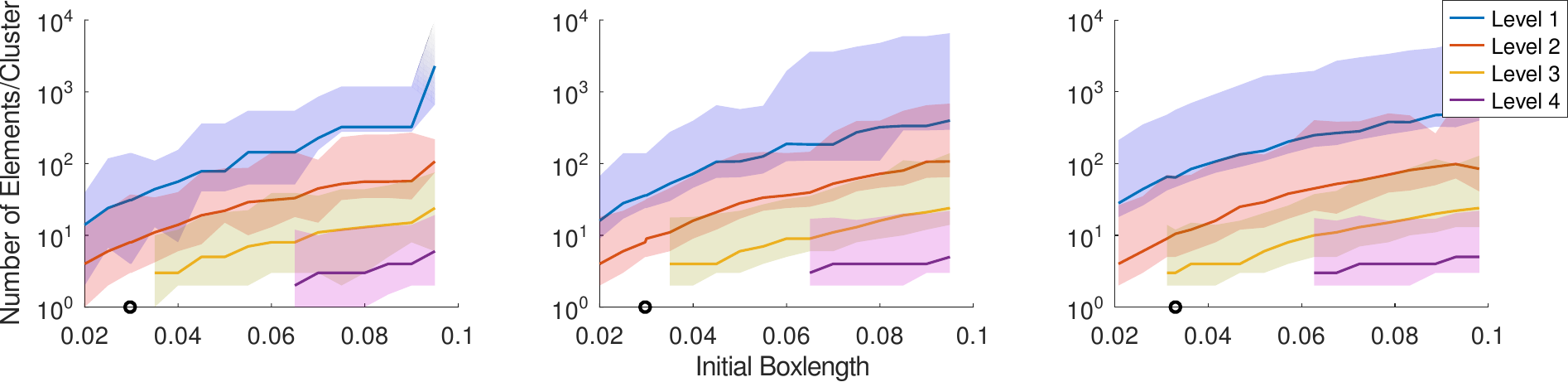}
  \caption{Median number of elements per cluster in each level plus lower and upper 10\% quantiles (shaded area) for the \emph{coarse} head mesh and the three clustering methods: Box1, Box2, and k-means. The circles represents the initial box length for the default setting. }\label{Fig:LukasNCcoarse}
\end{figure}

%%%%%%%%%%%%%%%%%%%%%%%%%%%%%%%%%%%%%%%%%%
\begin{table}[!h]
  \begin{center}
  \begin{tabular}{lccccc||c}\\
    & Min & 1st Quant & Median & 3rd Quant & Max & Mean\\
    Root Box1 &2 & 18.25 & 31.0 & 43 &7726 &152.96 \\ 
    Root Box2 &18 & 28 & 36.0 & 53 &7585 &195.28 \\ 
    Root k-means &38 & 51 & 70.0 & 348.75 &767 &195.28 \\ 
    \hline
    Leaf Box1 &1 & 6 & 8.0 & 12 & 1752 & 42.36 \\ 
    Leaf Box2 &3 & 6 & 8.0 & 12 & 1819 & 46.34 \\ 
    Leaf k-means &1 & 3 & 3.0 & 5 & 767 & 17.05 \\  
  \end{tabular}
  \end{center}
  \caption{Minimum, 1st Quantile, Median, 3rd Quantile, maximum, and mean value for the number of elements in the clusters in the root (upper part) and leaf levels for the \emph{coarse} head mesh.}\label{Tab:Coarse05}
\end{table}
When looking at the number of elements per cluster (see Fig.~\ref{Fig:LukasNCcoarse}) it becomes clear that for the coarse mesh a huge variation can be observed for all three clustering methods which is also reflected in the difference between mean and median number of elements per cluster, see Table \ref{Tab:Coarse05}. This difference  becomes even bigger in the leaf level that contains clusters that are small in radius but contain a high number of elements.

\begin{figure}[!h]
  \includegraphics[width=0.98\textwidth]{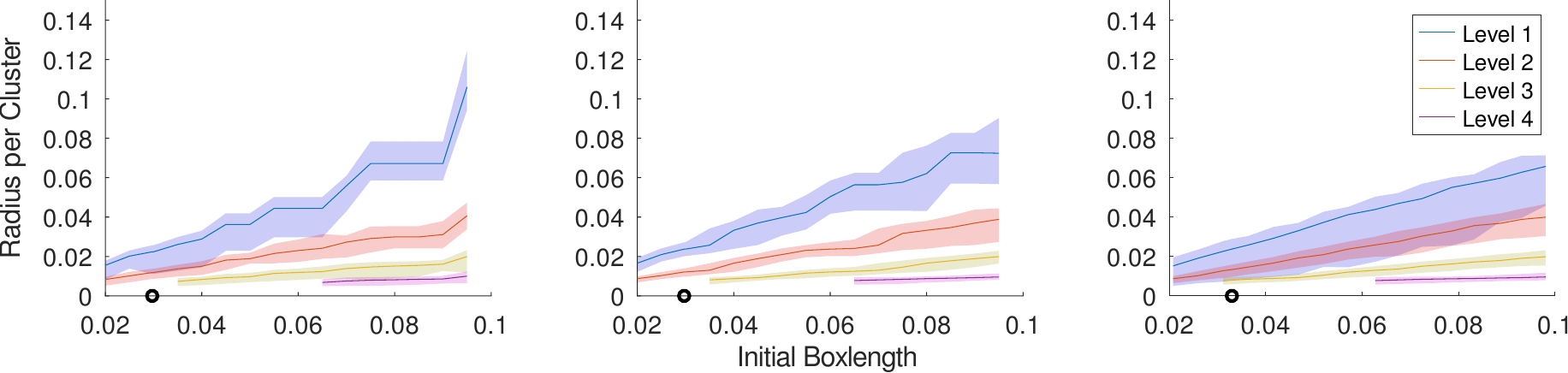}
  \caption{Median cluster radii in each level plus lower and upper 10\% quantiles (shaded area) for the \emph{coarse} head mesh and the three clustering methods: Box1, Box2, and k-means. The circles represents the initial box length for the default setting.}\label{Fig:Lukasrcoarse}
\end{figure}

Compared to the number of elements per cluster the cluster radii of all three clustering methods are rather balanced, see Fig.~\ref{Fig:Lukasrcoarse}. For the k-means based clustering there is a bigger variation in the root level which can be explained by the property of the k-means algorithm that was already observed for the spherical meshes. 

\begin{figure}[!h]
  \includegraphics[width=0.98\textwidth]{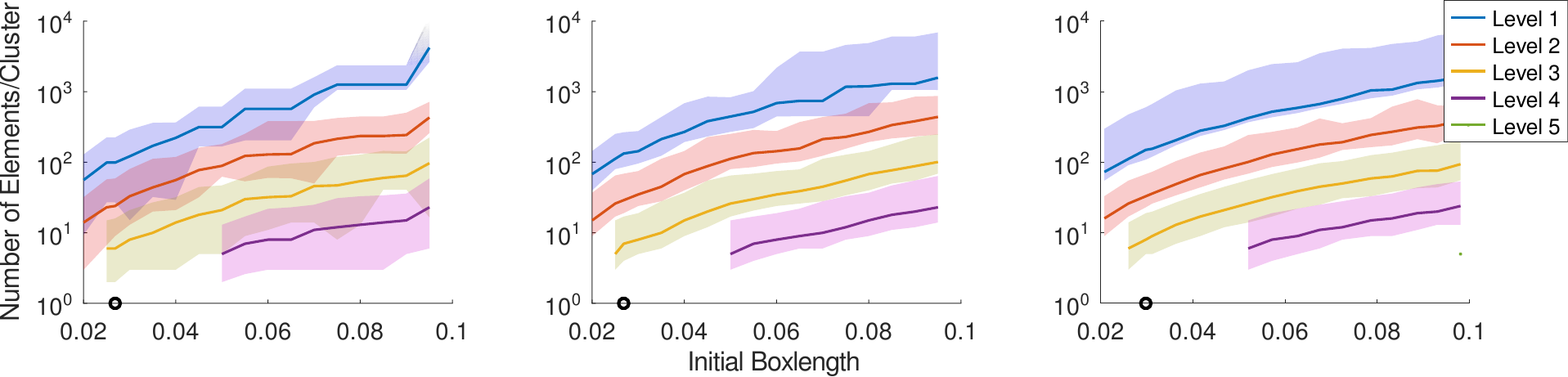}
  \caption{Median number of elements per cluster in each level plus lower and upper 10\% quantiles (shaded area) for the \emph{fine} head mesh and the three clustering methods: Box1, Box2, and k-means. The circles represents the initial box length for the default setting.}\label{Fig:LukasNCfine}
\end{figure}

\begin{figure}[!h]
  \includegraphics[width=0.98\textwidth]{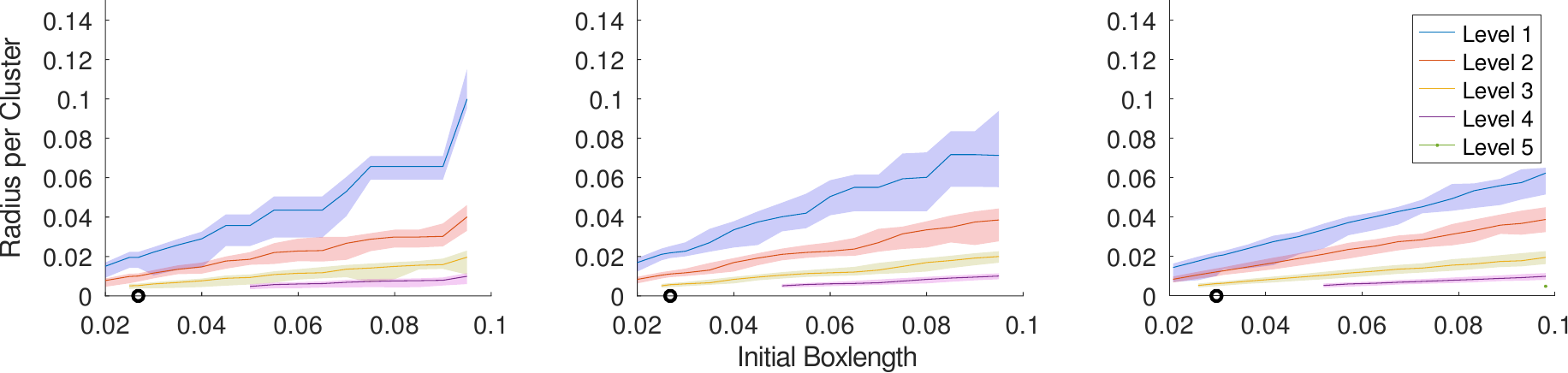}
  \caption{Median cluster radii in each level plus lower and upper 10\% quantiles (shaded area) for the \emph{fine} head mesh and the three clustering methods: Box1, Box2, and k-means. The circles represents the initial box length for the default setting.}\label{Fig:Lukasrfine}
\end{figure}

For the fine mesh (see Figs.~\ref{Fig:LukasNCfine} and \ref{Fig:Lukasrfine}, and Table \ref{Tab:Fine05}) the number of elements per cluster are much more balanced. One big difference in the clustering for the two meshes (see Figs.~\ref{Fig:LukasNCcoarse} to \ref{Fig:Lukasrfine}) is that for the default box lengths the finer mesh already has 3 multipole levels for all clustering methods. For the coarse setting the number of levels at the default box length for the box type clustering methods  is still 2, in case of the k-means clustering the default setting yields 3 levels. As already observed for the spherical meshes the additional level has positive influence on the performance of the MLFMM. 
\begin{table}[!h]
  \begin{center}
  \begin{tabular}{lccccc||c}\\
    & Min & 1st Quant & Median & 3rd Quant & Max & Mean\\
    Root Box1 & 15 & 54 & 99.0 & 128.75 &4350 &186.13 \\ 
    Root Box2 & 60 & 102 & 133.0 & 179 &6122 &262.46 \\ 
    Root k-means & 102 & 127 & 141.0 & 200.75 &986 &239.58 \\ 
    \hline
    Root Box1 & 1 & 4 & 6.0 & 8 & 447 & 12.41 \\ 
    Root Box2 & 1 & 5 & 7.0 & 10 & 510 & 14.58 \\ 
    Root k-means & 2 & 6 & 8.0 & 11 & 499 & 16.94 \\ 
  \end{tabular}
  \end{center}
  \caption{Minimum, 1st Quantile, Median, 3rd Quantile, maximum, and mean value for the number of elements in the clusters in the root (upper part) and leaf levels for the \emph{fine} mesh.}\label{Tab:Fine05}
\end{table}

%The high variation in the number of elements per cluster is partly because of the meshing itself and partly because of the geometry. For example around the ear lobe a cluster may contain twice as many elements as for cluster containing a smooth part of the mesh. 
%
%
%
\begin{figure}[!h]
  \begin{tabular}{rcrc}
    \raisebox{0.300\textwidth}{\small a) \hspace{-16pt}} &
    \includegraphics[width=0.44\textwidth]{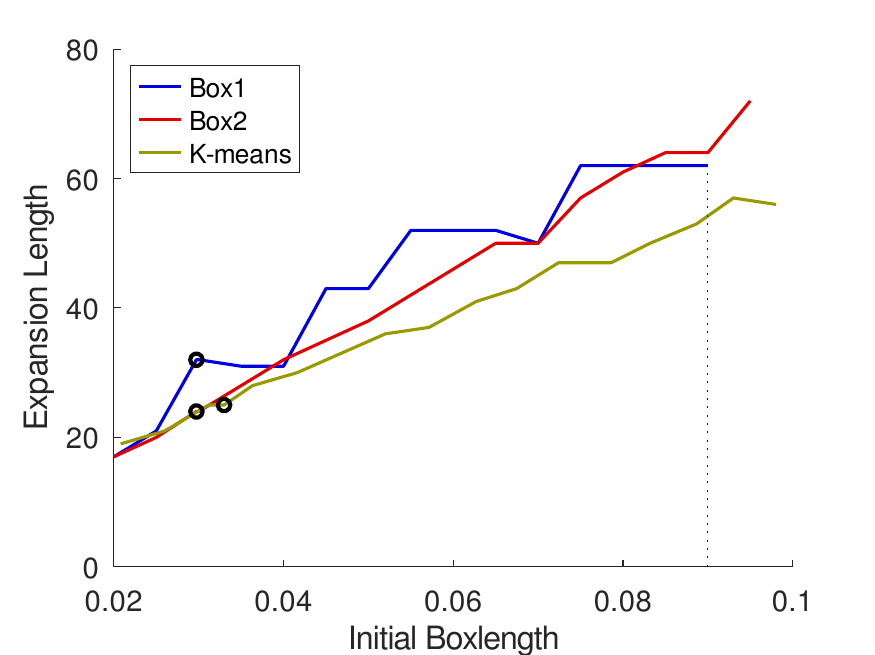} &
    \raisebox{0.300\textwidth}{\small b) \hspace{-16pt}} &
    \includegraphics[width=0.44\textwidth]{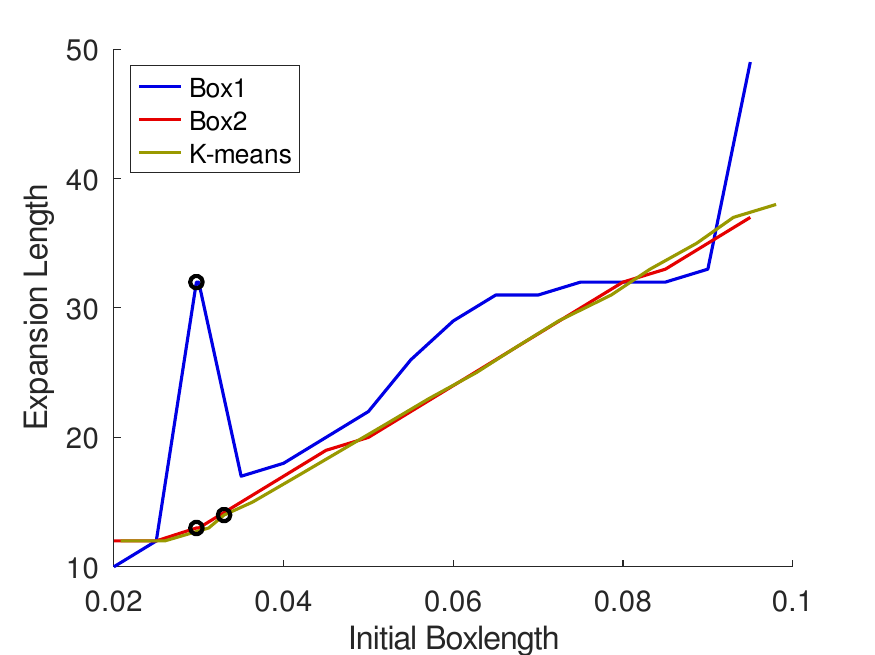}\\
    \raisebox{0.300\textwidth}{\small c) \hspace{-16pt}} &
    \includegraphics[width=0.44\textwidth]{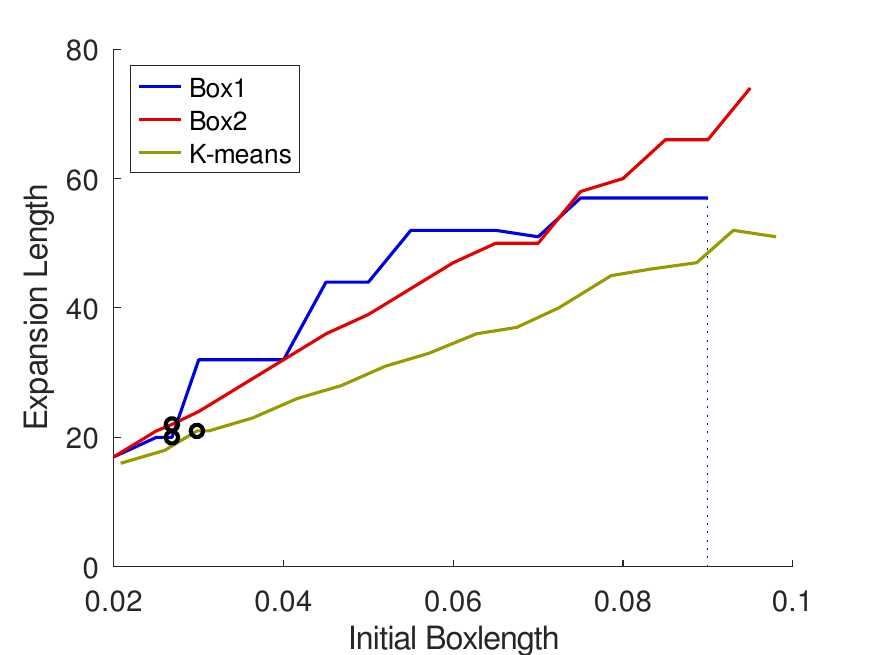} &
    \raisebox{0.300\textwidth}{\small d) \hspace{-16pt}} &
    \includegraphics[width=0.44\textwidth]{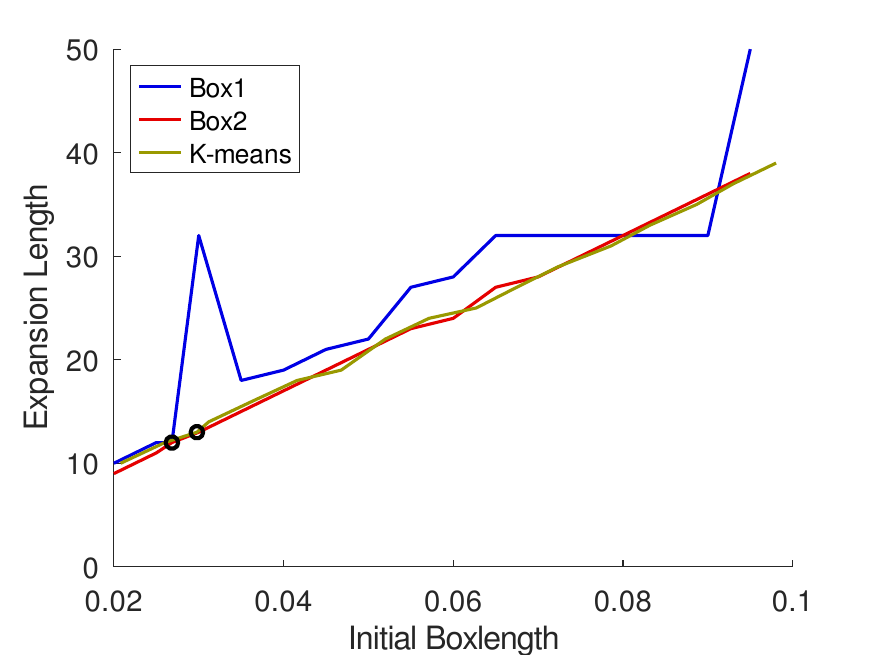}
  \end{tabular}
  \caption{Expansion length for the coarse (first row) and fine head meshes (second row) at the root level (first column) and the second multipole level (second column) at $f = 20000$\,Hz.}\label{Fig:LukasExplength}
\end{figure}

When looking at the expansion lengths (see Fig.~\ref{Fig:LukasExplength} for both meshes) an anomaly around $h_0 = 0.03$\,m for the Box1 clustering can be observed at the curves for level 2. The expansion lengths at the higher levels are the same as the expansion length at the root level. This anomaly is caused by a part of the pinna that is close to the ear canal and where the elements in the cluster at the root level are not contiguous (see Fig.\ref{Fig:Problems}b for a schematic representation of such a case). For the next level this parent cluster was split into two parts, however, since one part was only containing one element, both parts were merged again resulting in a single cluster with large radius and thus a large expansion length at the child level. This large expansion length causes several numerical problems, and the iterative solver for the BEM system did not converge. Note, that although this problem just occurs for the Box1 clustering in this case, there may be still a very small probability that this problem can also affect the other two clustering methods. However, the construction using a fixed subdivision in the root level for the Box1 clustering raises the probability of this problem to some degree.

\begin{table}
  \begin{center}
    \begin{tabular}{lccccc||c}\\
      & Min  &1st Quant & Median & 3rd Quant & Max & Mean\\
      Root Box1 &2 & 18.25 & 31.0 & 43 &7726 &152.96 \\ 
      Root Box2 &18 & 28 & 36.0 & 53 &7585 &195.28 \\ 
      Root k-means &38 & 51 & 70.0 & 348.75 &767 &195.28 \\ 
      \hline
      Root Box1 &1 & 1 & 2.0 & 3 & 493 & 11.24 \\ 
      Root Box2 &1 & 3 & 3.0 & 4 & 505 & 16.02 \\ 
      Root k-means &1 & 3 & 3.0 & 5 & 767 & 17.05 \\
    \end{tabular}
  \end{center}
  \caption{ Minimum, 1st Quantile, Median, 3rd Quantile, maximum, and mean value for the number of elements in the clusters in the root (upper part) and leaf levels for the coarse mesh with modified number of levels.}\label{Tab:Coarse03}
\end{table}

\subsubsection{Performance}\label{Sec:PerformanceHead}
\begin{figure}[!h]
  \includegraphics[width=0.46\textwidth]{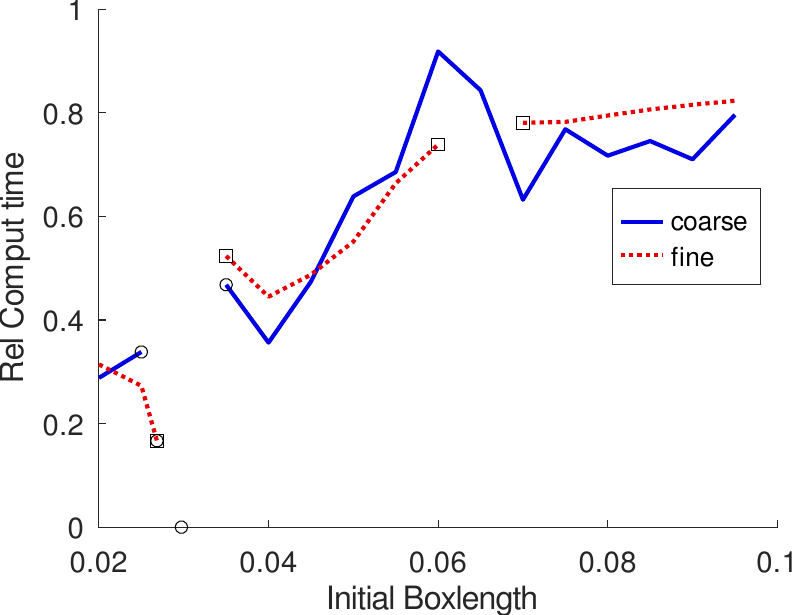}
  \includegraphics[width=0.48\textwidth]{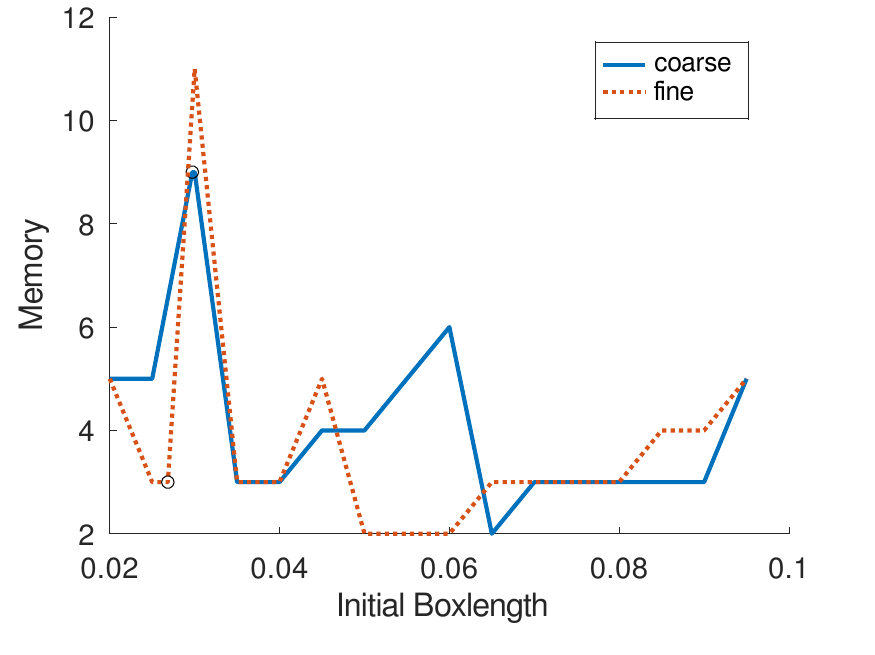}
  \caption{Relative computation time and peak memory in Gb for the two head meshes using the Box1 clustering. The dotted lines represent the values for the fine mesh.}\label{Fig:LukasTimeIter}
\end{figure}

In Fig.~\ref{Fig:LukasTimeIter} the differences in time and memory consumption for the two meshes and the Box1 clustering method at frequencies $f = 5000, 10000, 15000,$ and $20000$\,Hz are compared.
For both meshes the lines depicting the relative computation time contain a hole around $h_0\approx 0.03$\,m where the iterative solver for the global system did not converge. The reason for these instabilities was the one single non-contiguous element that raised the expansion length for the higher levels. There is also a gap for the fine mesh at $h_0 = 0.065$\,m, that cannot be attributed to a sudden rise in expansion length. This problem is caused by the mixture of small and large elements and by the bias of the  cluster midpoint towards a group of small elements, see also Fig.~\ref{Fig:Problems}c. This leads to few big outliers with relatively large cluster radii and an expansion length that is too big for the rest of the cluster interactions. %Fig~\ref{Fig:Lukasboxproblem}.
%    \begin{figure}[!h]
%      \includegraphics[width=0.8\textwidth]{Pics/ExpLengthLukas}
%      \caption{Expansion lengths for the two head meshes, where the dotted lines represent the expansion lengths for the finer mesh. The 'o' markers indicate the clusterlengths where a level was added for the coarse mesh, the 'x' marker denotes the equivalent values for the fine mesh.}
%    \end{figure}
The far-field criterion for two ``normal'' clusters is too small for the high order of the Hankel function, see also Section~\ref{Sec:Problems}.

In this specific case the numerical problem can be avoided by calculating the cluster midpoint as an area weighted sum over all vertices in the cluster, which reduces  the expansion length from 52 to 41.  An alternative way to deal with this problem is to adapt the expansion length to the cluster-pair involved in the FMM and not the largest cluster in the given level.
%
%  \begin{figure}[!h]
%    \includegraphics[width=0.5\textwidth]{Pics/Lukasboxproblem}
%    \caption{Boxplot of the clusterradii for $h_0 = 0.065$ and $f = 20000$ Hz if the cluster midpoint is calculated as the sum over all vertices in the cluster compared to the case when the midpoint is calculated as a weighted sum.}\label{Fig:Lukasboxproblem}
%  \end{figure}
%
%

\begin{figure}[!h]
  \begin{tabular}{rcrc}
    \raisebox{0.320\textwidth}{\small a) \hspace{-16pt}} &
    \includegraphics[width=0.44\textwidth]{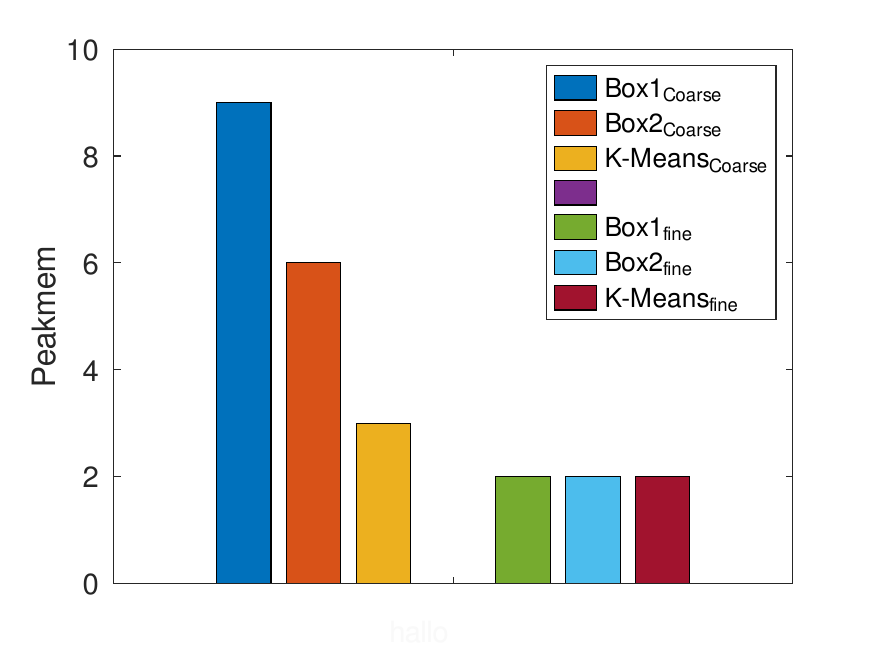} &
    \raisebox{0.320\textwidth}{\small b) \hspace{-16pt}} &
    \includegraphics[width=0.44\textwidth]{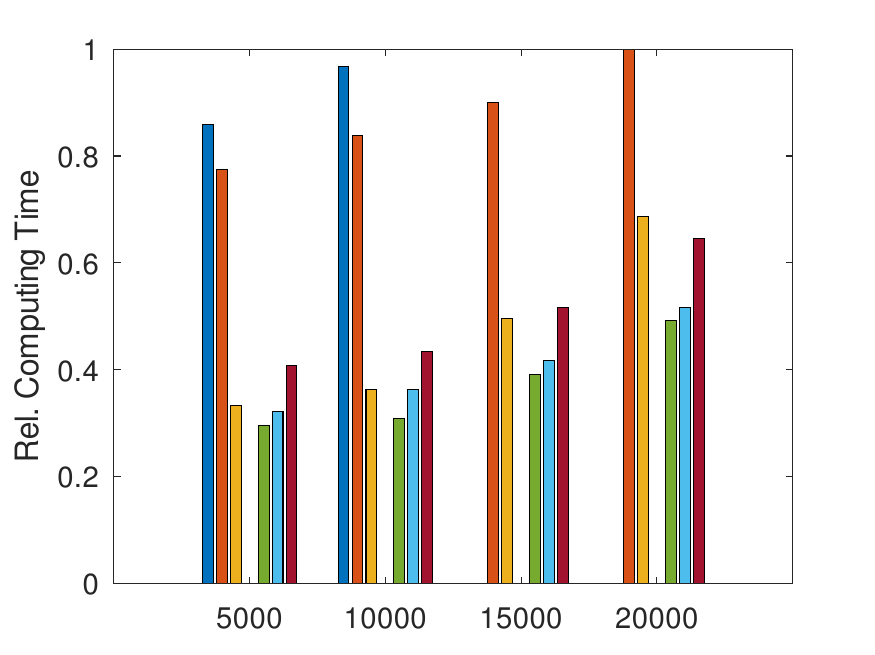} \\
    \raisebox{0.320\textwidth}{\small c) \hspace{-16pt}} &
    \includegraphics[width=0.44\textwidth]{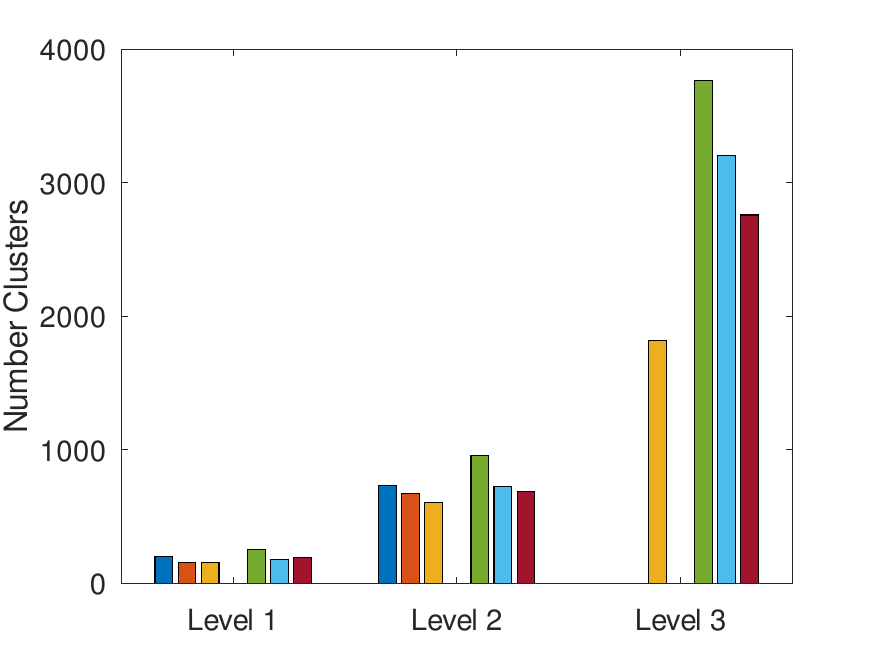} &
    \raisebox{0.320\textwidth}{\small d) \hspace{-16pt}} &
    \includegraphics[width=0.44\textwidth]{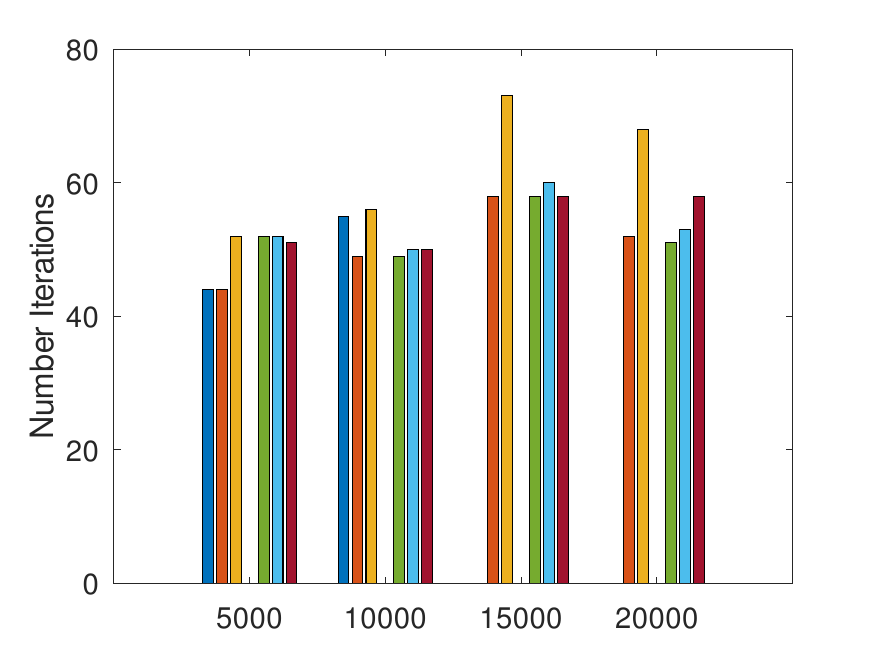}
  \end{tabular}
  \caption{a) Peak memory for the coarse and fine head meshes for the three clustering approaches, b) relative computation time for different frequencies, c) number of cluster per level, and d) number of iterations for each frequency for the three different clustering approaches. For all calculations the default setting for $h_0$ was used.}\label{Fig:CompHeadAll}
\end{figure}

In Fig.~\ref{Fig:CompHeadAll} the three clustering methods for the two head meshes are directly compared to each other. It becomes apparent that although having more elements the fine mesh, that is more balanced terms of element sizes, performs better with respect to memory consumption and computation time. For the fine mesh, the memory consumption does not differ very much between the clustering methods, while for the coarse mesh the number of levels and thus the number of elements in the leaf clusters has an influence on the performance. With the additional level for the k-means clustering method, the memory and computer time is much smaller then for the rest of the clustering approaches.
Note that for a better comparison the number of CPUs and threads has been set to 1. 
As already mentioned numerical problems arose for the Box1 clustering at 15000 and 20000 Hz. Computation time and number of iterations were not displayed for these types of clustering at these frequencies. Figs.~\ref{Fig:CompHeadAll}b) and d) indicate that although the number of iterations for $f = 5000$ Hz is almost the same for the different clustering methods, the additional level provides a much smaller matrix-vector multiplication time for the k-means clustering approach. For the fine mesh the Box based clustering methods resulted in 3 level cluster trees for all frequencies and the performance is also slightly better compared with the k-means based approach. Note, that in case of the fine mesh, the default box length was ``luckily'' small enough in the fine case not to fall into the region where instabilities occur, see Fig. \ref{Fig:LukasTimeIter}, where the default box length is depicted with a circle around $h_0 = 0.023$\,m.

If the rounding process in the approach to determine the levels (see Eq.~(\ref{Equ:Nlevels})) is slightly modified by rounding up the number of levels if the value after the decimal point is bigger or equal to 0.3 instead of 0.5, the Box1 and Box2 clustering methods for the head meshes become efficient again. The number of levels is raised to three, and the clusters at leaf level are more balanced (see Table  ~\ref{Tab:Coarse03}) although the stability problems for the Box1 clustering still persist.

%As Tab.~\ref{Tab:Coarse03} shows, the clusters in the leaf level are very unbalanced in terms of number of elements per clusters.
%
\begin{figure}[!h]
  \includegraphics[width=0.46\textwidth]{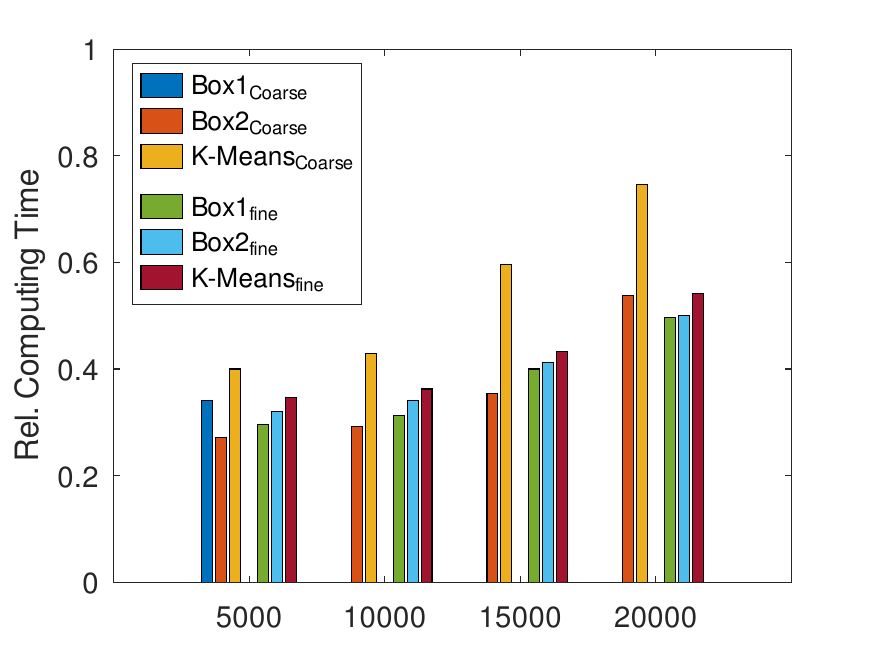}
  \includegraphics[width=0.46\textwidth]{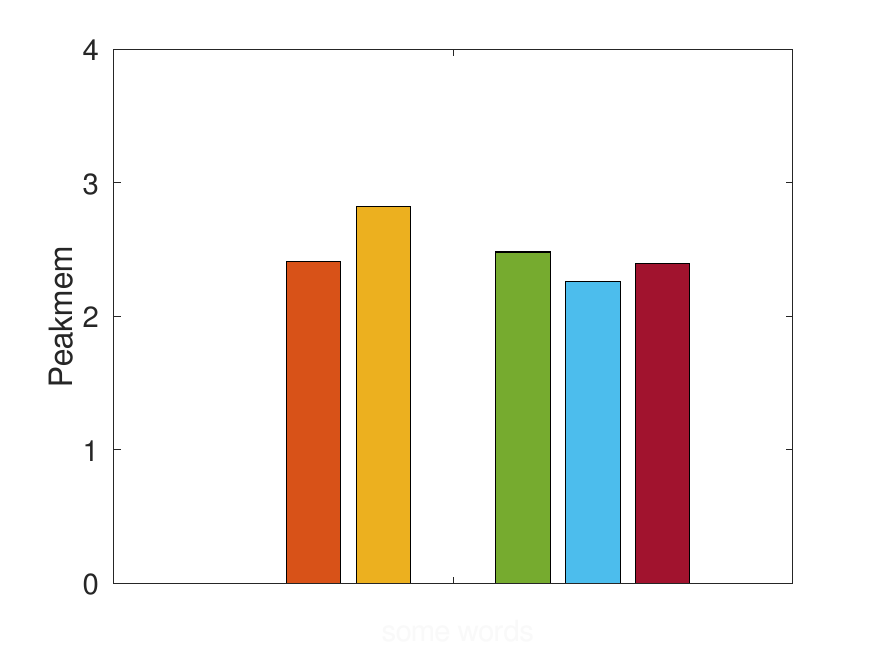}
  \caption{Relative computation time and peak memory for different clustering methods and the head meshes for the modified number of levels.}
\end{figure}
In summary, it is apparent that a lower of number of elements does \emph{not} always enhance the efficiency of the MLFMM. This is mostly because the clustering is based on average element sizes, but with non-uniform meshes the clustering results in leaf clusters containing many elements. This in turn, effects the number of entries in the nearfield matrix and has negative effect on the efficiency of the method.
%
%
%
%
%%%%%%%%%%%%%%%%%%%%%%%%%%%%%%%%%%%%%%%%%%%%%%%%%%%%%%%%%%%%%%%
\section{Discussion}\label{Sec:Conclusion}
In this paper different strategies for defining clusters for the Helmholtz MLFMM have been presented and compared using numerical experiments with 5 different meshes. Table \ref{Tab:Meshsummary} shows a comparison of these meshes  in terms of number of elements and their edge lengths.
\begin{table}[!h]
  \begin{center}
  \begin{tabular}{l | >{$}c<{$} | >{$}c<{$} | >{$}c<{$} | >{$}c<{$} }
    Name & \# \text{ Edges} & \text{Mean} & \text{Var} & \text{Var/Mean} \\
    \hline
    Uniform Sphere & 104430 & 0.02 & 2.16\cdot10^{-6} & 1.1\cdot10^{-4}\\
    Weak Non-Uni & 109839 & 0.018 & 4.65\cdot10^{-5} & 2.5\cdot 10^{-3}\\
    Strong Non-Uni & 104844 & 0.013 & 2.34\cdot10^{-4} & 1.7\cdot10^{-2}\\
    \hline
    Head coarse & 46575 & 2.12\cdot10^{-3} & 6.13\cdot10^{-6} & 2.8\cdot10^{-3}\\
    Head fine & 70077 & 2.23\cdot10^{-3} & 2.12\cdot10^{-6} & 9.5\cdot10^{-4}
  \end{tabular}
  \end{center}
  \caption{Number of edges, mean, variance and variance/mean of the edge lengths for the 5 meshes used.}\label{Tab:Meshsummary}
\end{table}

In a first numerical experiment the dependence of the MLFMM on the size of the clusters at the root level of the cluster tree has been investigated for uniform as well as non-uniform meshes of the unit sphere. Independent of the actual clustering method used having about $O(\sqrt{N})$ clusters at the root level is a good choice for building up the root of the cluster-tree.
This choice was first introduced for the single level fast multipole method \cite{Coifmanetal93}) and also proofs to be a good choice for the multi level approach.

When looking for the optimal clustering there is always the issue of the balance between efficiency and stability. On the one hand, having many clusters at root level implies a bigger computational effort, because for example the cluster to cluster interactions at the root level scale with $O(n_0^2 L^2) \approx O(n_0^2 N/n_0) = O( N n_0)$. But on the other hand, there is a gain in stability because more clusters mean smaller cluster radii and thus smaller expansion lengths $L$. It is important to point out, that there is \emph{not} a one value fits all solution to the problem of finding an optimal size for the root clusters, there is a big dependence on the uniformity of the mesh. 

This dependency also became apparent in the second round of numerical experiments that were made using the real life example of calculating the sound pressure at a human head and pinna. For this type of non-uniform meshes two different settings were investigated. The first mesh (coarse head mesh) was representation of a human head with $N = 31050$ elements with very small elements in the area of the pinnae but relatively large elements at the rest of the head. For the second (fine) mesh with $N = 46718$ elements the large elements on the head were subdivided once to provide a more balanced mesh in terms of element sizes. For the coarse mesh the variance of the cluster radii at leaf level is relatively moderate, but the number of elements per cluster  ranges from 1 to 1752 elements with a mean of 42.36, while the median is much lower (see Tab.~\ref{Tab:Coarse05}). This explains the high effort of the MLFMM for non-uniform meshes. The (few) clusters with a large number of elements negatively influence the effort for setting up and solving the costly nearfield part. 

Non-uniform meshes may also suffer from stability problems because the truncation order of the MLFMM and in turn the maximum order of the Hankel functions necessary for the multipole expansion depends on the largest radius of all clusters at a given level. The argument of the Hankel functions on the other hand depends on the distance and local radii of the two clusters for which the FMM is calculated. If a level contains large and small clusters this fact may lead to numerical problems this may lead to numerical problems, because the argument of the Hankel function is too small for the high order.

An option to avoid such a problem would be to adapt the far-field criterion Eq.~(\ref{Equ:FarfieldCrit}) to use the biggest cluster radius on the level instead of individual local radii, but this would result in unnecessarily high computational effort. Alternatively, one can adapt the FMM procedure to adapt the truncation of the FMM process (i.e., the expansion length) to the cluster radii, while still keeping the quadrature nodes for the sphere the same for all clusters to avoid interpolating nodes on the sphere on a single level. This approach was successfully applied to the head meshes to avoid the convergence problems.

Also, as illustrated by the comparison between the coarse and the fine head mesh, the number of elements in a mesh is not necessarily the main criterion for the efficiency of the BEM calculations when the multilevel FMM is involved, but rather how the cluster tree is built up. Even with the modified number of levels the times needed for the computations for the coarse and fine head meshes are in the same range, and there is almost no gain in choosing a mesh with fewer elements over a balanced mesh.
\section*{Acknowledgments}
We would like to thank Felix Perfler and Lukas Thalhammer for providing the original head mesh used for the HRTF calculations. This work was supported by the European Union (EU) within the project SONICOM (grant number: 101017743, RIA action of Horizon 2020).
\bibliography{Mesh2hrtf}

@article{Ziegelwangeretal16,
title = {A priori mesh grading for the numerical calculation of the head-related transfer functions},
optlongjournal = {Applied Acoustics},
journal = {Appl. Acoust.},
volume = {114},
pages = {99--110},
year = {2016},
issn = {0003-682X},
doi = {https://doi.org/10.1016/j.apacoust.2016.07.005},
opturl = {https://www.sciencedirect.com/science/article/pii/S0003682X1630192X},
optlongauthor = {Harald Ziegelwanger and Wolfgang Kreuzer and Piotr Majdak},
author = {H. Ziegelwanger and W. Kreuzer and P. Majdak},
optkeywords = {Head-related transfer functions, Boundary element method, Mesh grading},
optabstract = {Head-related transfer functions (HRTFs) describe the directional filtering of the incoming sound caused by the morphology of a listener’s head and pinnae. When an accurate model of a listener’s morphology exists, HRTFs can be calculated numerically with the boundary element method (BEM). However, the general recommendation to model the head and pinnae with at least six elements per wavelength renders the BEM as a time-consuming procedure when calculating HRTFs for the full audible frequency range. In this study, a mesh preprocessing algorithm is proposed, viz., a priori mesh grading, which reduces the computational costs in the HRTF calculation process significantly. The mesh grading algorithm deliberately violates the recommendation of at least six elements per wavelength in certain regions of the head and pinnae and varies the size of elements gradually according to an a priori defined grading function. The evaluation of the algorithm involved HRTFs calculated for various geometric objects including meshes of three human listeners and various grading functions. The numerical accuracy and the predicted sound-localization performance of calculated HRTFs were analyzed. A-priori mesh grading appeared to be suitable for the numerical calculation of HRTFs in the full audible frequency range and outperformed uniform meshes in terms of numerical errors, perception based predictions of sound-localization performance, and computational costs.}
}

@Article{Coifmanetal93,
  author = 	 {R. Coifman and V. Rokhlin and S Wandzura},
  title = 	 {The Fast Multipole Method for the Wave Equation: A Pedestrian Prescription},
  optlongjournal = 	 {IEEE Antennas and Propagation Magazine},
  journal =      {IEEE Trans. Antennas Propag.},
  year = 	 {1993},
  OPTkey = 	 {},
  volume = 	 {35},
  number = 	 {3},
  pages = 	 {7--12},
  OPTmonth = 	 {},
  OPTnote = 	 {},
  OPTannote = 	 {}
}

@Article{Rahola96,
  optlongauthor =       {Jussi Rahola},
  author = {J. Rahola},
  title =        {Diagonal forms of the translation operators in the
                 fast multipole algorithm for scattering problems},
  journal =      {BIT},
  volume =       {36},
  number =       {2},
  pages =        {333--358},
  optmonth =        jun,
  year =         "1996",
  optCODEN =        "BITTEL, NBITAB",
  optDOI =          "http://www.springerlink.com/openurl.asp?genre=article&id=doi:10.1007/BF01731987",
  optISSN =         "0006-3835 (print), 1572-9125 (electronic)",
  optMRclass =      "65R20 (65F10 78-08 78A45)",
  optMRnumber =     "97k:65298",
  optMRreviewer =   "Wojciech Mydlarczyk",
  optbibdate =      "Wed Jan 4 18:52:24 MST 2006",
  optbibsource =    "ftp://ftp.math.utah.edu/pub/bibnet/subjects/fastmultipole.bib;
                 ftp://ftp.math.utah.edu/pub/tex/bib/bit.bib;
                 http://springerlink.metapress.com/openurl.asp?genre=issue&issn=0006-3835&volume=36&issue=2",
  optURL =          "http://www.mai.liu.se/BIT/contents/bit36.html;
                 http://www.springerlink.com/openurl.asp?genre=article&issn=0006-3835&volume=36&issue=2&spage=333",
  optacknowledgement = ack-nhfb,
  optkeywords =     "fast multipole method",
}

@article{CecDar13,
optlongauthor = {Cecka, Cris and Darve, Eric},
author = {Cecka, C. and Darve, E.},
title = {Fourier Based Fast Multipole Method for the {H}elmholtz Equation},
optlongjournal = {SIAM Journal on Scientific Computing},
journal = {SIAM J. Sci. Comput.},
volume = {35},
number = {1},
pages = {A79-A103},
year = {2013},
optdoi = {10.1137/11085774X},
optURL = {https://doi.org/10.1137/11085774X},
opteprint = {https://doi.org/10.1137/11085774X}
}

@article{Rokhlin90,
 author = {Rokhlin, V.},
 title = {Rapid solution of integral equations of scattering theory in two dimensions},
 journal = {J. Comput. Phys.},
 volume = {86},
 number = {2},
 year = {1990},
 optissn = {0021-9991},
 pages = {414--439},
 optdoi = {http://dx.doi.org/10.1016/0021-9991(90)90107-C},
 optpublisher = {Academic Press Professional, Inc.},
 optaddress = {San Diego CA, USA},
 }

@article {Rokhlin93,
	author = "V. Rokhlin",
	title = "Diagonal Forms of Translation Operators for the {H}elmholtz Equation in Three Dimensions",
	optlongjournal = "Applied and Computational Harmonic Analysis",
	journal = "Appl. Comput. Harmon. Anal.",
	volume = "1",
	year = "December 1993",
	optabstract = "<P>The diagonal forms are constructed for the translation operators for the Helmholtz equation in three dimensions. While the operators themselves have a fairly complicated structure (described somewhat incompletely by the classical addition theorems for the Bessel functions), their diagonal forms turn out to be quite simple. These diagonal forms are realized as generalized integrals, possess straightforward physical interpretations, and admit stable numerical implementation. This paper uses the obtained analytical apparatus to construct an algorithm for the rapid application to arbitrary vectors of matrices resulting from the discretization of integral equations of the potential theory for the Helmholtz equation in three dimensions. It is an extension to the three-dimensional case of the results of Rokhlin (<B>J. Complexity</B> <B><I>4</I></B>(1988), 12-32), where a similar apparatus is developed in the two-dimensional case.<B>Copyright 1993, 1999 Academic Press</B></P>",
	pages = "82-93",
	opturl = "http://www.ingentaconnect.com/content/ap/ha/1993/00000001/00000001/art01006"
}

@misc{mesh2hrtf,
  title = {{Mesh2HRTF -- An open source software for the numerical calculation of head-related transfer functions}},
  author = {{The Mesh2HRTF Developers}},
  note = {(accessed 2022-12-20)},
  url = {https://mesh2hrtf.org}
}

@inproceedings{Ziegelwangeretal15,
optlongauthor = {Ziegelwanger, Harald and Kreuzer, Wolfgang and Majdak, Piotr},
author = {Ziegelwanger, H. and Kreuzer, W. and Majdak, P.},
year = {2015},
OPTmonth = {07},
pages = {},
title = {{Mesh2HRTF}: An open-source software package for the numerical calculation of head-related transfer functions},
booktitle = {Proceeding of the 22nd International Congress on Sound and Vibration},
address = {Florence, Italy},
doi = {10.13140/RG.2.1.1707.1128}
}

@book{SauSch10,
  title={Boundary Element Methods},
  author={Sauter, S.A. and Schwab, C.},
  isbn={9783540680932},
  series={Springer Series in Computational Mathematics},
  year={2010},
  publisher={Springer Berlin Heidelberg}
}

@Article{Brinkmannetal23,
  title = {Recent {{Advances}} in an {{Open Software}} for {{Numerical HRTF Calculation}}},
  author = {Brinkmann, F. and Kreuzer, W. and Thomsen, J. and Dombrovskis, S. and Pollack, K. and Weinzierl, S. and Majdak, P.},
  year = {2023},
OPTmonth = jul,
  volume = {71},
  number = {7/8},
  pages = {502--514},
  doi = {10.17743/jaes.2022.0078},
  urldate = {2023-07-18},
  journal = {J. Audio Eng. Soc.}
}

@Article{AmiPro00,
  author = 	 {S. Amini and A. Profit},
  title = 	 {Analysis of the truncation errors in the fast multipole method for scattering problems},
  journal = 	 {J. Comput. Appl. Math.},
  year = 	 {2000},
  OPTkey = 	 {},
  volume = 	 {115},
  OPTnumber = 	 {},
  pages = 	 {23--33},
  OPTmonth = 	 {},
  OPTnote = 	 {},
  OPTannote = 	 {}
}

@Article{AmiPro03,
  author = 	 {S. Amini and A.T.J. Profit},
  title = 	 {Multi-level fast multipole solution of the scattering problem},
  journal = 	 {Eng. Anal. Bound. Elem.},
  year = 	 {2003},
  OPTkey = 	 {},
  volume = 	 {27},
  OPTnumber = 	 {},
  pages = 	 {547--564},
  OPTmonth = 	 {},
  OPTnote = 	 {},
  OPTannote = 	 {}
}

@Article{Darve00,
  author = 	 {E. Darve},
  title = 	 {The Fast Multipole Method: Numerical Implementation},
  journal = 	 {J. Comput. Phys.},
  year = 	 {2000},
  OPTkey = 	 {},
  volume = 	 {160},
  OPTnumber = 	 {},
  pages = 	 {195--240},
  OPTmonth = 	 {},
  OPTnote = 	 {},
  OPTannote = 	 {},
  doi =          {10.1006/jcph.2000.6451}
}

@article{Lietal24,
title = {Novel and efficient implementation of multi-level fast multipole indirect BEM for industrial Helmholtz problems},
optjournal = {Engineering Analysis with Boundary Elements},
journal = {Engn. Anal. Bound. Elem.},
volume = {159},
pages = {150--163},
year = {2024},
issn = {0955-7997},
doi = {10.1016/j.enganabound.2023.11.027},
opturl = {https://www.sciencedirect.com/science/article/pii/S0955799723005568},
author = {Y. Li and O. Atak and W. Desmet},
}

@article{Kreuzeretal24,
title = {NumCalc: An open-source {BEM} code for solving acoustic scattering problems},
optlongjournal = {Engineering Analysis with Boundary Elements},
journal = {Engn. Anal. Bound. Elem.},
volume = {161},
pages = {157--178},
year = {2024},
issn = {0955-7997},
doi = {10.1016/j.enganabound.2024.01.008},
opturl = {https://www.sciencedirect.com/science/article/pii/S0955799724000171},
optauthorlong = {Wolfgang Kreuzer and Katharina Pollack and Fabian Brinkmann and Piotr Majdak},
author = {W. Kreuzer and K. Pollack and F. Brinkmann and P. Majdak},
}

@Article{Greengardetal98,
  author = 	 {L. Greengard and J. Huang and V. Rokhlin and S. Wandzura},
  title = 	 {Accelerating Fast Multipole Methods for the Helmholtz Equation at Low Frequencies},
  journal = 	 {IEEE Comput. Sci. Eng.},
  year = 	 {1998},
  OPTkey = 	 {},
  volume = 	 {5},
  number = 	 {3},
  pages = 	 {32--38},
  OPTmonth = 	 {},
  OPTnote = 	 {},
  OPTannote = 	 {}
}

@manual{octave,
    title     = {{GNU Octave} version 9.4.0 manual: a high-level interactive language for numerical computations},
    author    = {J.W. Eaton and D. Bateman and S. Hauberg and R. Wehbring},
    year      = {2024},
    url       = {https://www.gnu.org/software/octave/doc/v9.4.0/},
  }

@Article{Chengetal06,
  author = 	 {H. Cheng and W.Y. Crutchfield and Z. Gimbutas and L.F. Greengard and J. F. Ethridge and J. Huang and V. Rokhlin and N. Yarvin and J. Zhao},
  title = 	 {A wideband fast multipole method for the {H}elmholtz equation in three dimensions},
  journal = 	 {J. Comput. Phys.},
  year = 	 {2006},
  OPTkey = 	 {},
  volume = 	 {216},
  OPTnumber = 	 {},
  pages = 	 {300-325},
  OPTmonth = 	 {},
  OPTnote = 	 {},
  OPTannote = 	 {}
}

@Article{Yunetal20,
  author = 	 {D. Yun and H. Jung and H. Kang and D. Seo},
  title = 	 {{Acceleration of the Multi-Level Fast Multipole Algortihm Using K-Means Clustering}},
  journal = 	 {Electronics},
  year = 	 {2020},
  OPTkey = 	 {},
  volume = 	 {9},
  OPTnumber = 	 {},
  OPTpages = 	 {},
  OPTmonth = 	 {},
  OPTnote = 	 {},
  OPTannote = 	 {}
}

@Article{Moeller92,
  author = 	 {H. M\o{}ller},
  title = 	 {Fundamentals of binaural technology},
  journal = 	 {Appl. Acoust.},
  year = 	 {1992},
  OPTkey = 	 {},
  volume = 	 {36},
  OPTnumber = 	 {},
  pages = 	 {171--218},
  OPTmonth = 	 {},
  OPTnote = 	 {},
  OPTannote = 	 {}
}

@article{WangYu26,
title = {A statistical shape model–based framework for virtual head generation and head-related transfer function database construction},
journal = {Applied Acoustics},
volume = {254},
pages = {111400},
year = {2026},
issn = {0003-682X},
doi = {https://doi.org/10.1016/j.apacoust.2026.111400},
OPTurl = {https://www.sciencedirect.com/science/article/pii/S0003682X26001817},
optauthor = {Yewei Wang and Guangzheng Yu},
author = {Y. Wang and G. Yu},
OPTabstract = {Personalized head-related transfer function (HRTF) databases are essential for virtual auditory display (VAD) applications, yet conventional databases rely on human measurements or detailed anatomical scans that involve personal biometric information and associated ethical concerns. Virtual-head-based HRTF databases therefore offer a scalable and privacy-preserving alternative. This study proposes a statistical shape model (SSM)-based framework for generating anatomically plausible virtual heads that conform to the anatomical distribution of an existing population dataset and for computing their corresponding HRTFs via numerical simulation. High-resolution three-dimensional scans of 100 Chinese adults were acquired and processed to construct a fixed-topology representation while preserving acoustically relevant morphological variability, particularly in the pinna regions. Principal component analysis (PCA) was applied to construct the SSM and characterize inter-subject shape variations. By sampling modal coefficients from a multivariate Gaussian distribution in the retained SSM space and applying an automatic 99}
}
\bibliographystyle{elsarticle-num}
\end{document}